\documentclass[12pt]{article}
\usepackage{amssymb,amsmath,amsfonts,mathrsfs,upgreek} 
\usepackage{mathabx}

\usepackage{graphicx} 
\usepackage{bm}

\usepackage[colorlinks=true]{hyperref}
\usepackage{comment}

\renewcommand{\marginpar}[2][]{} 

\renewcommand\comment[1]{{\iffalse #1 \fi}}

\newtheorem{theorem}{Theorem}[section]

\newtheorem{lemma}{Lemma}[section]
\newtheorem{proposition}{Proposition}[section]

\newtheorem{remark}{Remark}[section]
\newtheorem{notation}{Notation}[section]

\newcommand{\io}{{\infty}}

\newcommand{\real}{ {\mathbb R}   }
\newcommand{\torus}{ {\mathbb T}    }
\newcommand{\integer}{ {\mathbb Z}   }
\newcommand{\complex}{ {\mathbb C}   }
\newcommand{\bB}{ {\mathbb B}   }

\newcommand{\cB}{ {\mathcal B}   }

\newcommand{\cP}{ {\mathcal H}   }
\newcommand{\cM}{ {\mathcal M}   }

\renewcommand{\Im}{\, {\rm Im}\,}
\renewcommand{\Re}{\, {\rm Re}\,}

\newcommand{\eproof}{\qed}

\newcommand\beq[1]{ \begin{equation}\label{#1} }
\newcommand{\eeq}{ \end{equation} }
\newcommand{\beqno}{ \[ }
\newcommand{\eeqno}{ \] }
\newcommand\beqa[1]{ \begin{eqnarray} \label{#1}}
\newcommand{\eeqa}{ \end{eqnarray} }
\newcommand{\beqano}{ \begin{eqnarray*} }
\newcommand{\eeqano}{ \end{eqnarray*} }

\newtheorem{definition}{Definition}[section]
\newcommand\dfn[1]{ \begin{definition}\label{#1} }
\newcommand\edfn{ \end{definition} }
\newcommand\notat[1]{ \begin{notation} \label{#1} 
 }
\newcommand\enotat{\end{notation}}

\newcommand\rem{\begin{remark} 
\rm 
}
\newcommand\erem{\end{remark} 
}

\newcommand{\proof}{\par\medskip\noindent{\bf Proof\ }}
\newcommand\equ[1]{{\rm (\ref{#1})}}
\newcommand{\nl}{{\smallskip\noindent}}
\newcommand{\giu}{{\medskip\noindent}}
\newcommand{\Giu}{{\bigskip\noindent}}
\newcommand{\noi}{{\noindent}}
\newcommand{\qed}{\hskip.5truecm
\vrule width 1.7truemm height 3.5truemm depth 0.truemm
\par\Giu}
\newcommand{\qedeq}{\hskip.5truecm
\vrule width 1.7truemm height 3.5truemm depth 0.truemm}

 \newcommand\casialt[3]{ \left\{  \begin{array}{ll}
 {#1} & \mbox{ {\rm if} ${#2}$} \\
 {#3} & \mbox{ {\rm otherwise}}
 \end{array} \right.}

\newcommand{\e}{\varepsilon}

\renewcommand{\a }{\alpha }

\newcommand{\s }{\sigma }
\newcommand{\ii }{{\rm i} }
\renewcommand{\d }{\delta }

\newcommand{\g }{\gamma}
\newcommand{\f }{\varphi}

\renewcommand{\L }{\Lambda }
\newcommand{\m }{\mu }

\renewcommand{\t }{\tau }
\renewcommand{\o }{\omega }
\renewcommand{\O }{\Omega }
\newcommand{\C}{\mathbb{C}}
\newcommand{\Z}{\mathbb{Z}}

\newcommand{\normadue}{|}

\newcommand{\cgot}{\mathfrak c}

\def\N{\mathbb N}
\def\R{\mathbb R}
\def\T{\mathbb T}

\def\const{{\, \rm const\, }}
\def\dst{\displaystyle}
\def\bks{\, \backslash\, }
\def\meas{{\rm\, meas\, }}
\def\Tp{{T_K^\perp}}
\def\TNp{T_N^\perp}
\newcommand\eqby[1]{\stackrel{\equ{#1}}{=}}

\newcommand\proiezione{\, {\mathtt p}}

\newcommand\modulo{|}
\newcommand\ham{H_\e}
\newcommand\sa{\theta} 
\newcommand\loge{\lambda}

\newcommand{\tetta}{\vartheta }

\newcommand{\K}{K_{{}_2}}
\newcommand{\KO}{K_{{}_1}}
\newcommand\ks[1]{\kappa_\varstar(#1)}

\newcommand{\bs}{{\bar s}}
\newcommand{\ts}{{\tilde s}}

\newcommand{\norma}{\thickvert\!\!\thickvert}
\usepackage{appendix}

\title{\bf 
On the topology of nearly--integrable Hamiltonians at simple resonances
\ \\
}

\begin{document}

\author{ 
\footnotesize L. Biasco  \& L. Chierchia
\\ \footnotesize Dipartimento di Matematica e Fisica
\\ \footnotesize Universit\`a degli Studi Roma Tre
\\ \footnotesize Largo San L. Murialdo 1 - 00146 Roma, Italy
\\ {\footnotesize biasco@mat.uniroma3.it, luigi@mat.uniroma3.it}
\\ 
}

\maketitle

\begin{abstract}

\noi
We show that, in general,  averaging  at simple resonances  a  real--analytic, nearly--integrable Hamiltonian, one obtains a one--dimensional system with a cosine--like potential; ``in general'' means for a generic class of holomorphic perturbations and apart from a finite number of simple resonances with small Fourier modes; 
 ``cosine--like'' means that the potential depends only on the resonant angle, with respect to which it is a Morse function with one  maximum and one  minimum.
\\
Furthermore, the (full) transformed Hamiltonian  is the sum
of  an effective one--dimen\-sio\-nal Hamiltonian (which is, in turn,  the sum of the unperturbed Hamiltonian  plus the  cosine--like potential) and  a perturbation, which is exponentially small with respect to the oscillation of the potential. 
\\
As a corollary,  under the above hypotheses, if the unperturbed Hamiltonian is also strictly convex, the effective Hamiltonian at {\sl any simple resonance} (apart a finite number of low--mode resonances) has the phase portrait of a pendulum.   
\\
The results presented in this paper  are an essential step in the proof (in the ``mechanical'' case)
 of a conjecture by 
Arnold--Kozlov--Neishdadt (\cite[Remark~6.8, p. 285]{AKN}), 
claiming that the  measure of the ``non--torus set'' in general nearly--integrable Hamiltonian systems
has the same
 size of the perturbation; compare \cite{BClin}, \cite{BC}.

\end{abstract}

{\small 
\tableofcontents
}




\newcommand\gen{{\cal G}^n_1}
\def\genk{{\cal G}^n_{1,K}}
\def\genKO{{\cal G}^n_{1,\KO}}
\def\genK{{\cal G}^n_{1,\K}}

\renewcommand{\cgot}{\gamma}
\newcommand{\npiu}{{\bar n}}
\newcommand\noruno[1]{  |#1|_{{}_1} }
\newcommand\norinf[1]{  |#1|_{{}_\infty} }
\renewcommand\ln{\log}
\newcommand\st{{\rm \ s.t.\ }}
\newcommand\tail{{\tau}}
\renewcommand\ks[1]{\tail_{\rm o}(#1)}
\newcommand{\ells}{{\ell^\infty(\integer^n_\varstar)}} 
\newcommand\hol{{\bB}}
\newcommand\palla{{\bf B_{\rm 1}}}

\newcommand\Lk{{\rm S}_k}
\newcommand\Ak{{\rm A}_k}
\newcommand\hAk{{\rm \hat A}_k}
\newcommand\pk{{\proiezione^\perp_k}}
\def\dist{{\rm\, dist\, }}
\newcommand\Dk{D^{1,k}}
\newcommand\tDk{{\tilde D}^{1,k}}
\newcommand\tailo{\tail_{\rm o}}
\newcommand\gk{\f_k}

\section{Introduction}

Consider a real--analytic, nearly--integrable Hamiltonian given, in action--angle variables, by
\beq{ham}
H_\e(y,x)=h(y)+\e f(y,x)\ ,\quad  (y,x)\in \cM:=D\times \torus^n \ ,
\eeq
where $D$ is a bounded domain in $\real^n$, $\torus^n=\real^n/(2\pi\integer^n)$ is the usual flat $n$ dimensional torus and $\e$ is a small parameter measuring the size of the perturbation $\e f$. The phase space  $\cM$ is endowed with the standard symplectic form $dy\wedge dx$ so that the Hamiltonian flow  $\phi_{H_\e}^t(y_0,x_0)=:(y(t),x(t))$ governed by $H_\e$ is the solution of the standard Hamiltonian equations
\beq{prosciutto}
\left\{\begin{array}l
\dot y = -\partial_x H_\e(y,x)\ ,\\
\dot x= \partial_y H_\e(y,x)\ ,
\end{array}
\right.\qquad 
\left\{\begin{array}l
y(0)=y_0\ ,\\
x(0)=x_0\ ,
\end{array}
\right.
\eeq
(where $t$ is time and dot is time derivative).

\nl
It is well known that, in general, the $\phi_{H_\e}^t$--dynamics  is strongly influenced by 
{\sl resonances of the (unperturbed) frequencies} 
$\o(y):=h'(y)=\partial_y h(y)$, 
i.e., by rational relations $$\o(y)\cdot k=\sum_{j=1}^n \o_j(y)k_j=0\,,$$
with $k\in \integer^n\bks\{0\}$; for general information, compare, e.g.,  \cite{AKN}.
Indeed, assuming a standard KAM non--degeneracy assumption on $h$, e.g., that the frequency map  $y\in D\to \o(y)$
 is a real--analytic diffeomorphism of $D$ onto the ``frequency space'' $\O:=h'(D)$, then
the action space $D$ can   be covered by three open sets 
\beq{coperta}
D\subseteq D^0\cup D^1\cup D^2
\eeq
so that the following holds. Roughly speaking, $D^0\times\torus^n$ is a fully non--resonant set which is filled, up to an exponentially  small set,
by {\sl primary KAM tori}, namely, by homotopically trivial, Lagrangian tori $\phi^t_{H_\e}$--invariant  on which the flow is  analytically conjugated to the linear flow  
$$\theta\in\torus^n\mapsto \theta+\o t\ $$
with $\o$ satisfying a Diophantine 
condition
\beq{diofanto}
|\o\cdot k|\ge\frac{ \gamma}{\noruno{k}^\tau} \,\qquad  \forall \ k\in \integer^n\bks\{0\} \ ,
\eeq
(for some $\gamma,\tau>0$); $\o\cdot k$ denoting the standard inner product $\sum \o_i k_i$ and 
$\noruno{k}:=\sum |k_i|.$ Furthermore, such tori are deformation of integrable tori.
 \\
$D^1$ is an open $O(\sqrt\e)$--neighbourhood of {\sl simple resonances} (i.e., of regions where {\sl exactly one} independent resonance $\o(y)\cdot k=0$ holds) and 
$D^2$ is a set of measure $O(\e)$; compare the Covering Lemma  (Proposition~\ref{beccalossi}) below\footnote{This description follows by choosing carefully certain parameters (such as the ``small divisor constant'' $\a$  and  ``Fourier cut--offs'' $K$) as functions of $\e$ and  disregarding logarithmic corrections.}.

\nl
The region $D^2$ contains double (and higher) resonances and, in general, in 
$D^2\times\torus^n$ there are $O(\e)$ regions where the dynamics  is  {\sl non--perturbative}, being  ``essentially'' governed  (after suitable rescalings)  by  an $\e$--independent Hamiltonian; compare\footnote{\label{toronto}
\cite[Remark 6.8, p. 285]{AKN}: ``{\sl It is natural to expect that in a generic system with three or more degrees of freedom the measure of the ``non--torus'' set has order $\e$. Indeed, the $O(\sqrt\e)$--neighbourhoods of two resonant surfaces intersect in a
domain of measure $\sim \e$. In this domain, after the partial averaging taking into account the resonances under consideration, normalizing the deviations of the ``actions'' from the resonant values by the quantity $\sqrt\e$, normalizing time, and discarding the terms of higher order, we obtain a Hamiltonian of the form $1/2(Ap, p) + V (q_1, q_2)$, which does not involve a small parameter (see the definition of the quantity $p$ above). Generally speaking, for this Hamiltonian there is a set of measure $\sim 1$ that does not contain points of invariant tori. Returning to the original variables we obtain a ``non--torus'' set of measure $\sim \e$.}'' }
\cite{AKN}.

\nl
The dynamics in the simple--resonance region $D^1\times\torus^n$ is particularly relevant and interesting. For example, it  plays a major r\^ole in {\sl Arnold diffusion}, as showed by Arnold himself \cite{A},  who based his famous instability argument on shadowing partially hyperbolic trajectories arising near simple resonances. \\
On the other hand,  in $D^1\times\torus^n$ there appear {\sl secondary KAM tori}, namely $n$--dimensional  KAM tori with different topologies, which depend upon specific characteristics of the perturbation $\e f$. The appearance of  secondary tori is a genuine non--integrable effect, since such tori do not exist in the integrable regime.\\
 In the announcement \cite{BClin}  it is claimed that,  in the case of mechanical systems -- namely systems governed by Hamiltonians of the form $|y|^2/2+\e f(x)$ --  and for generic potentials $f$,  
primary and       
 secondary tori fill the region $D^1\times\torus^n$ up to a set of measure nearly exponentially small, showing that the ``non--torus set'' is, at most, $O(\e)$ as conjectured in\footnote{See footnote~\ref{toronto} above.
 Note also that, as it was proved in \cite{L}
(in dimension 2) and
\cite{Nei},
\cite{poschel1982}   (in any dimension),  the union of  {\sl primary}
invariant tori fills the phase space up to a set of measure   $O(\sqrt{\e})$.
This result is optimal: the phase region inside the separatrix of the   pendulum  
$\frac12 y^2 +\e \cos x$, with $y\in \real$ and $x\in\torus^1$,
   does not contain any 
primary invariant torus, namely a circle which is a global graph over the angle on $\torus^1$, and this region has measure $4\sqrt{2\e}$.
Indeed this region is filled by secondary tori, corresponding to oscillations of the pendulum.
 } \cite{AKN} and studied in \cite{MNT}.
 In fact, Theorem~\ref{scozia} below is one of the building block of the proof (in the mechanical case) 
 of the Arnold--Kozlov--Neishdadt conjecture as outlined in \cite{BC}.

\giu
This paper is devoted to the fine topological and quantitative analysis of  the behaviour of generic  systems in the simple resonant region $D^1\times\torus^n$.

\giu
In this introduction, we  briefly  discuss  the main aspects of this analysis in the particular case of {\sl purely positional potentials}; precise statements are given in Theorem~\ref{scozia} of \S~\ref{mainresult} below, and, for the general (but more technical and  implicit) case, in Theorem~\ref{scoziaY} of \S~\ref{denteazzurro}).

\nl
$D^1$ is the union of suitable regions $\Dk$, which  are $O(\sqrt \e)$-close to {\sl exact simple resonances} $\{y\in D|\  \o(y)\cdot k=0\}$, and which are labelled by {\sl generators} $k$
of one dimensional, maximal sublattices of $\integer^n$ (see \eqref{cippalippa} below); ``exact'' meaning that 
$\o(y)$ does not verify  double or higher resonant relations. \\
In averaging (or normal form) theory, one typically considers a {\sl finite } but large (possibly, $\e$--dependent) number of simple resonances. More precisely, one considers generators $k$ with 
$\noruno{k}\le K$,  and $K$  can be chosen according to the application one has in mind. Typically, one chooses $K\sim 1/\e^a$ for a suitable $a>0$
(as in Nekhoroshev theorem \cite{poschel}, \cite{GCB}) or 
$K\sim |\log \e|^a$ (as in the KAM theory for secondary tori of \cite{BClin}, \cite{BC}).

\nl
By averaging theory, in any fixed simple resonant region $\Dk$, one can remove the non--resonant angle dependence, so  as
to  symplectically conjugate $H_\e$, for $\e$ small enough, 
to a  Hamiltonian of the form
\beq{varano}
  h(y)+\e G^k(y,k\cdot x) +\e  R^k(y,x)
\eeq
where $\sa\to G^k(y,\sa)$ is a function of one angle 
and
$R^k$ is a ``very small'' remainder
with\footnote{$R^k_\ell(y)$ denotes the $\ell^{\rm th}$ Fourier coefficient
of $x\mapsto R^k(y,x).$} $R^k_{jk}(y)=0, \forall j\in\integer$.
Thus, up to the remainder $R^k$, the Hamiltonian depends effectively only on
the ``resonant angle''
$\sa:=k\cdot x$ and therefore the ``effective Hamiltonian'' $h+\e G^k$ is {\sl integrable}: this is the starting point for (``a priori stable'') Arnold diffusion or for the KAM theory for secondary tori of \cite{BClin}--\cite{BC}.

\nl
Obviously, there are here two main issues: 
\begin{itemize}
\item[\bf (a)] {\sl What is the actual ``generic form'' of $G^k$}?

\item[\bf (b)]  {\sl How small (and compared to what) is the remainder $R^k$}?
\end{itemize}

\nl
{\bf (a)} According to averaging (or normal form)  theory\footnote{For generalities on Averaging Theory, see, e.g., \cite[\S~6]{AKN} and references therein.} $G^k$ is ``close''  to the 
projection of the potential $f$ on the Fourier modes of the resonant maximal sublattice $k\Z$:
\beq{proiezione}
\proiezione_{k\Z}f(x)=\sum_{j\in\integer } f_{jk} e^{\ii jk\cdot x}\,.
\eeq
Now, since $f$ is real-analytic on $\torus^n$, it is holomorphic in a complex strip $\torus^n_s$ around $\torus^n$ of width $s>0$ and its Fourier coefficients decay exponentially fast as\footnote{Precise norms will be introduced in the next section \S~\ref{mainresult}.} $|f_k|\sim \|f\| e^{-|k|s}$. Hence, {\sl typically} (i.e., if $f_k\neq 0$) 
$$
\sum_{j\in\integer } f_{jk} e^{\ii jk\cdot x} = 
f_k e^{\ii k\cdot x} + f_{-k} e^{-\ii k\cdot x} + O \big(\|f\| e ^{-2|k| s}\big)
$$
which, by the reality condition $f_{-k}=\bar f_k$, can be written as
\beq{piedini}
\sum_{j\in\integer } f_{jk} e^{\ii jk\cdot x} = 
2|f_k|
 \cos\big(k\cdot x +\sa^{(k)}\big) + O \big(\|f\| e ^{-2|k| s}\big)
\eeq
for a suitable $\sa^{(k)}\in [0,2\pi)$. Thus,
\beq{genericus}
\proiezione_{k\Z}f(x)= \sum_{j\in\integer } f_{jk} e^{\ii jk\cdot x} = 
2|f_k|
\Big(  \cos\big(k\cdot x +\sa^{(k)}\big) +o(1)\Big)
\eeq
{\sl provided} 
\beq{pavone}
|k| \gtrapprox \frac1{s} \ .
\eeq
In other words,  one expects \equ{genericus} to hold for {\sl generic real--analytic potentials} and for all generators $k$'s satisfying \equ{pavone}. \\
Indeed, this is the case: we shall introduce  certain  classes of periodic holomorphic functions $\cP_{s,\tail}$, for which (choosing suitably the ``tail'' fuction $\t$) \equ{genericus} holds for generators $k$ satisfying \equ{pavone}. \\
The class $\cP_{s,\tail}$ turns out to be  ``generic'' in several ways:\\ 
(i)  it contains an open dense set in the class of real--analytic functions having holomorphic extension on a complex neighbourhood  of size $s$ of $\torus^n$ (in the topology induced by a suitably weighted Fourier norm);
\\
(ii) its  unit ball  is of measure 1 (with respect to a  natural probability measure);
\\
(iii) it is a ``prevalent set''. 
\\
For precise statements see  Definition~\ref{orcadi} and Proposition~\ref{grizman} below.

\nl
Next,   in order for $G^k$ to be  close to $\proiezione_{k\Z}f$ in \equ{genericus}, one needs to have a bound of the type
\beq{rana}
\sup_{\Dk\times\torus^n}|
G^k-\proiezione_{k\Z}f|
\ll |f_k| \sim \|f\| e^{-\noruno{k}s}\,.
\eeq
As well known, averaging methods involve an {\sl analyticity loss} in  complex domains. In particular, 
 the Hamiltonian in \eqref{varano}
and, therefore, $G^k$,  can be analytically defined only in a {\sl smaller} complex
strip $\torus^n_{s_{\varstar}}$ with $s_{\varstar}<s$. Therefore, by analyticity arguments, the 
best one can hope for is an estimate of the type
\beq{fava}
\sup_{\Dk\times\torus^n}|
G^k-\proiezione_{k\Z}f|
\le c \cdot  \|f\| e^{-\noruno{k}s_{\varstar}}\ , 
\eeq
for a suitable constant  $c$ that can be taken to be smaller than any prefixed positive number.
But then, for \equ{fava} and \equ{rana} to be compatible one sees that one must  ``essentially'' have 
$s_\varstar\sim s$ and that standard averaging theory is  not enough\footnote{Compare, e.g., \cite{poschel}, where 
$s_{\varstar}=s/6$. For a more detailed 
comparison with the averaging lemma of  \cite{poschel}, 
see also Remark \ref{cb500x}--(iv) below.
Compare also \cite{delshamsG} and \cite{giorgilli}.}. To overcome this problem, we provide   (Section~\ref{vongole})
a normal form lemma with small analyticity loss, ``small'' meaning that one can take  
\beq{figo}
s_{\varstar}=s(1-1/K)
\eeq
(compare, in particular,  \eqref{dadostar}). The value \equ{figo} 
 is compatible with \eqref{rana} for
$\noruno{k}\leq K$, showing that, indeed, generically, one has
$$
G^k(y, k\cdot x)= 2|f_k|
\Big(  \cos\big(k\cdot x +\sa^{(k)}\big) +o(1)\Big)
\ .
$$
In particular we prove that: {\sl  The ``effective Hamiltonians'' $h+\e G^k$, as $k$ vary ($c/s\le |k|\le K$),  have (up to a phase--shift) the same cosine--like form}  and, hence, the same  topological feature; compare, also,  Remark \ref{kuno}--(i) below.

\nl
Notice also that on low modes this last property, in general, does not holds,
as one immediately sees by considering  $k=e_1=(1,0,\ldots,0)$ and a potential $f$ such that 
$$
\proiezione_{e_1\Z} f(x):= \cos x_1+  \cos 2 x_1\,,
$$
which is a Morse function with two maxima and two minima in $\torus^1$.

\nl
{\bf (b)} What we just discussed gives also an indication for the question ``with respect to what $R^k$ has to be small''.
In fact, if, as expected,  \equ{genericus} is the leading behaviour,   one should have
\beq{marylin}
\|R^k\|\ll  |f_k|\ .
\eeq
But in order to perform averaging procedures, one  has, typically, control on  small divisors up to the truncation order $K$, so that 
the remainder  will contain high Fourier modes,  $|k|\sim K$, of the potential $f$. Such terms are  bounded  by $O(e^{-Ks})$, which are of the {\sl same size} of $\|G^k\|$, at least for $|k|\lesssim K$. 

\nl
To overcome this problem, we introduce in \S~\ref{geometria}, at difference with standard geometry of resonances (such as in \cite{nek}, \cite{poschel}, \cite{GCB}),  {\sl two} Fourier cut offs $\K\ge 3 \KO$ 
in such a way that on the simple resonant regions  $\Dk$ one has  non--resonance conditions for double and higher resonances up to  order $\K$, while $\KO$ is the maximum value of the size of the generators $k$ (i.e., $|k|\le \KO$). Therefore, we will  get an estimate of the remainder $R^k$ of the type
\beq{katarina}
\|R^k\|\le C\,  \K^a e^{-\K s/2}\le C'\,  |f_k|\ e^{-\K s/8}\ ,
\eeq
for suitable constants $C,C'>0$.
The final upshot is the {\sl complete normal form} 
\beq{dile}
\ham\circ\Psi_k
=:
h(y)+
2|f_k|\e
\Big( \cos(k\cdot x +\sa^{(k)})+
\tilde G^k(y,k\cdot x)+ \tilde R^k(y,x)
 \Big)
\eeq
with $\|\tilde G^k\|\ll 1$ and $\| \tilde R^k\|\le  C'\,  e^{-\K s/8}$; compare Theorem~\ref{scozia} below and, in particular, formula \equ{HkTE}.

\Giu
Summarizing:  {\sl for all $k$ large enough, $G^k$ is ``cosine--like''  i.e. a Morse function with one  maximum and one  minimum (compare Remark \ref{kuno}--(i) below) and 
  $R^k$ is exponentially small with respect to the oscillations  of $G^k$
  (see \equ{HkTE} and \equ{lothlorienTE} below).}

\noi
As a consequence we get that, {\sl if $h(y)$ is strictly convex, 
the effective Hamiltonian  has a phase portrait of a pendulum
(compare Remark~\ref{kuno}--(iv) below)}.

\section{Statements}
\label{mainresult}

\nl
Assume that $H_\e$ in \equ{ham}, for some $r,s>0$, 
admits
holomorphic extension on the complex domain 
$D_r\times\torus^n_s$, where 
$D_r\subseteq\C^n$  is the open complex  neighbourhood of $D$ formed by points  $z\in \C^n$ such that\footnote{We denote by $|\cdot|$ the usual Euclidean norm.} $|z-y|<r$, for some $y\in D$
and $\torus^n_s$
denotes the open complex neighbourhood of $\torus^n$ given by
\beq{toroseduto}
\torus^n_s:=\{x=(x_1,...,x_n)\in\complex^n:\ |\Im x_j|<s\}/(2\pi \integer^n)\ .
\eeq
The integrable hamiltonian $h$ is supposed to be ``KAM non--degenerate'' in the following sense.

\nl
{\bf Assumption A} {\sl Let $h$ be a real-analytic function 
\beq{h}
h: y\in D\subset \real^n \mapsto h(y)\in \real\ ,\quad (n\ge 2)\ , 
\eeq
where $D$ is a bounded domain of $\R^n$ and such that the {\sl frequency map} 
\beq{omega}
y\in D\mapsto \o(y):=\partial_yh(y)\in\O:=\o(D)\subseteq B_M(0)\subset \R^n \ ,\qquad  M:=\sup_D|\o(y)|\ ,
\eeq
{\sl is a  global diffeomorphism of $D$ onto $\O$} with Lipschitz constants given by 
\beq{Lip}
|y-y_0| \bar L^{-1}\le 
|\o(y)-\o(y_0)|\le L |y-y_0|\ , \quad (\forall\ y,y_0\in D)\ .
\eeq}
%
%
%

\noindent
Now, we describe the covering of frequency/action domain, which allows
to apply averaging theory (Proposition~\ref{pesce} below) to a perturbation of an integrable system with Hamiltonian $h$ at non--resonant  (modulus a lattice) zones\footnote{We use here the term ``zone'' in a loose way, {\sl not} in the technical meaning of Nekhoroshev's Theory; compare, e.g., \cite{GCB}.}. \\
Here, the main point is to find a suitable covering of simple resonances, which are the regions where the averaged Hamiltonian is 
integrable\footnote{In the sense that the, up to a small remainder,
the averaged Hamiltonian depends  on one angle.}, up to a small remainder. All other higher--order resonances are covered by one set, which is of small
measure: how small depending on the choice of the various parameters involved and it will vary according to the applications one has in mind.

\nl
Let $ \integer^n_\varstar$ denote the set of integer vectors $k\neq 0$ in $\integer^n$ such that the  first non--null  component is positive:
\beq{coscia}
 \integer^n_\varstar:=
 \big\{ k\in\integer^n:\ k\neq 0\ {\rm and} \ k_j>0\ {\rm where}\ j=\min\{i: k_i\neq 0\}\big\}\ ,
 \eeq
and denote by $\gen$  the {\sl generators of 1d maximal lattices}, namely, the set of  vectors $k\in  \integer^n_\varstar$ such that the greater common divisor (gcd)  of their components is 1:
\beq{cippalippa}
\gen:=\{k\in \integer^n_\varstar:\ {\rm gcd} (k_1,\ldots,k_n)=1\}\ .
\eeq
Then, {\sl the list of one--dimensional maximal lattices is given by the sets $\integer k$ with $k\in \gen$}.
\\
Given $K>0$ we set\footnote{$\noruno{k}:=\sum_{1\leq i\leq n}|k_i|.$}
\begin{equation}\label{elrond}
\genk:=\gen \cap \{\noruno{k}\leq K \}\,.
\end{equation}

\begin{proposition}\label{beccalossi}{\bf (Covering Lemma)}
Let  $h$ and $\o$ be as in Assumption~{\bf A} (\S~\ref{mainresult})   and fix $\K\ge \KO\ge 2$ and $\a>0$. Then, the domain $D$ can be covered by three sets $D^i\subseteq D$, 
\beq{Di} 
D=D^0\cup D^1\cup D^2\ ,
\eeq
 so that the following holds.

\nl
{\rm (i)}
$D^0$ is $(\a/2,\KO)$ completely non--resonant (i.e., non--resonant modulus $\{0\}$), namely, 
\begin{equation}\label{cipollotto}
y\in D^0 \quad\Longrightarrow\quad
|\o(y)\cdot k|\geq \a/2\,,\ \ \ \forall\ \, 0<|k|_{{}_1}\leq \KO\,.
\end{equation}
{\rm (ii)} $\dst D^1=\bigcup_{k\in\genKO} \Dk$, where,  
for each\footnote{Recall \eqref{elrond}.} 
 $k\in
 \genKO,$
 $\Dk$ is a neighbourhood of a simple resonance $\{y\in D: \o(y)\cdot k=0\}$, which is 
 $(2\a\K/|k|,\K)$ non--resonant modulo $\integer k$, namely, 
\begin{equation}\label{cipollotto2}
y\in \Dk \quad\Longrightarrow\quad
|\o(y)\cdot \ell |\ge 2\a\K /|k|\,,
\ \ \ \forall\ \, \ \ell\in \integer^n\ ,\ \ell\notin \Z k\ , \    \noruno{\ell}\leq \K\ .
\end{equation}
{\rm (iii)}  $D^2$  contains all the resonances of order two or more and has Lebesgue measure small with $\a^2$: more precisely, 
there exists a constant $c>0$
depending only on $n$ such that
\begin{equation}\label{teheran4} 
\meas (D^2) \le c \  \bar L^n M^{n-2}\ \a^2\  \K^{n+1} \KO^{n-1}
\ .
\end{equation}
\end{proposition}

\rem\label{noccetti}
(i) The neighbourhoods $\Dk$ of simple resonances $\{y: \o(y)\cdot k=0\}$  are explicitly defined as follows.
Denote  by
$\proiezione^\perp_k$
the orthogonal
projection  on the subspace perpendicular  to\footnote{Explicitly, $\proiezione^\perp_k \o:= \o- \frac{1}{|k|^2}(\o\cdot k) k$.} 
$k$ and, for $k\in\genKO$, define 
\beq{sonno}
\O^{1,k}:= 
\Big\{\o\in\real^n: |\o\cdot k|<\a,\  |\pk \o|<M, \ {\rm and}\   |\pk \o\cdot \ell|> \frac{3 \a \K}{|k|}  \ , \forall
\ell\in \genK\bks\Z k\Big\}
\eeq
Then, 
\beq{fischiettobis}
\Dk:=\{y\in D: \o(y)\in \O^{1,k}\}
\,.
\eeq
(ii) The domains $D^0,  D^2$ are explicitly defined 
in  \eqref{fischietto}, \eqref{neva0} and \eqref{wetton} below.

\nl
(iii) The simply resonant regions $\Dk$  in the  above Proposition are labelled by generators of 1--d maximal lattices $k\in \gen$ up to size $\noruno{k}\le \KO$, however,
{\sl the non--resonance condition \equ{cipollotto2} holds for integer vectors $\ell$ with  $\noruno{\ell}$ up to a (possibly) larger order $\K$}. 
This  improvement (with respect to having $\K=\KO$
as, e.g., in \cite{poschel})  is technical but important if one wants to have sharp control over the averaged Hamiltonian in a normal form near simple resonances;
in particular in order to obtain
\eqref{552bis} and \eqref{lothlorien}, which lead to \eqref{lothlorienTE}.

\nl
(iv)
The non--resonance relations \equ{cipollotto} and \equ{cipollotto2} allow to apply averaging theory and to remove the dependence upon the ``non--resonant angle variables''  up to exponential order; for precise statements, see   Theorem \ref{raffo} in \S~\ref{vento}.
\erem
We proceed, now, to describe the generic non--degeneracy assumption on periodic holomorphic functions, which will allow to state the main theorem (for the case of positional potentials). 

\giu
If $s>0$,  we denote by $\hol_s^n$   the Banach space of
{\sl real--analytic} functions on $\torus^n_s$ having zero average and finite  $\ell^\io$--Fourier norm:
\beq{hol}
\hol_s^n:=\Big\{f=\sum_{k\in\Z^n \atop k\neq 0} f_k e^{\ii k\cdot x} :\ \|f\|_s
:=\sup_{k\in \integer^n} |f_k| e^{|k|_{{}_1}s}<\io\Big\}\ .
\eeq
Note that $f\in \hol_s^n$ can be  uniquely written as:
\beq{arold}
f(x)= \sum_{k\in \gen}
\sum_{j\in \integer\backslash\{0\}} f_{jk}e^{\ii j k\cdot x}
\eeq 
For  functions\footnote{Not necessarily holomorphic in $y$.}
 $f:D_r\times \torus^n_s\to\complex$ we
 will also use the (stronger) norm\footnote{See Remark \ref{ozono} for details.}
\begin{equation}\label{canuto}
\norma f\norma_{D,r,s}=
\norma f\norma_{r,s}:=
\sup_{y\in D_r}\sum_{k\in\Z^n} 
|f_k(y)| e^{|k|_{{}_1}s}\,.
\end{equation}


\begin{definition}[Non degenerate potentials]\label{orcadi}
A  tail function $\tail$ is, by definition,  a 
non--increasing,  non--negative 
continuous function  
$$\tail:\d\in (0,1] \mapsto \tail(\d) \ge 0\ .
$$

\nl
Given  
$ s>0
$
 and a (possibly $s$--dependent) tail function $\tail$, we 
 define, for $\d\in (0,1]$,  $\cP_{s,\tail}(\d)$  as the set of functions in $ \hol_s^n$ 
such that,   for  any generator $k\in\gen$, the following 
holds\footnote{\label{befana2} One could substitute $n$ with every $\npiu>n/2$
in \equ{star}; compare Remark \ref{befana} below. The ``weight'' $ |k|_{{}_1}^{-n}$ is necessary in order to show that $\cP_{s,\tail}(\d)$ 
in \equ{amsterdam} has positive measure in a suitable probability space; compare Proposition~\ref{grizman}--(ii) below.}
\beq{star}
\mbox{if} \  
|k|_{{}_1}> \tail(\delta),
\ \mbox{then} \quad  
|f_k|\geq \d |k|_{{}_1}^{-n}\ e^{-|k|_{{}_1}s}\ ,
\eeq
The class $\cP_{s,\tail}$ is the union over $\d$  of the classes $\cP_{s,\tail}(\d)$:
\beq{amsterdam}\displaystyle\cP_{s,\tail}:=\bigcup_{0<\d\le 1} \cP_{s,\tail}(\d)\ .
\eeq
\end{definition}
The classes $\cP_{s,\tail}$ contain (if the tail is choosen properly) the non degenerate potentials for which
Theorem~\ref{scozia} below holds and, as mentioned in the Introduction, satisfy three main genericity properties\footnote{Such properties hold for {\sl any} tail $\t$, which  can be chosen differently according to the particular problem at hand.}, as showed in  
Proposition~\ref{grizman}  below  (compare, also, Remark \ref{pioggia}).

\begin{theorem}\label{scozia}
Let $n\ge 2$, $s>0,$ $0<\d,\gamma \le 1$  such that
\beq{porticato3}
\g\d<\frac{2^9}{s^{n}} e^{-{n}^2/2} 
\,.
\eeq
Consider a Hamilonian $H_\e(y,x)=h(y)+\e f(x)$ as in \eqref{ham} 
where $h$ satisfies
the non--degeneracy Assumption~{\bf A} (\S~\ref{mainresult}) and
$f$ is purely positional (i.e., independent of the $y$--variable) 
with
\begin{equation}\label{badaTE}
   \|f\|_{s}=1\,  .
\end{equation}
Assume  that the potential is   non--degenerate in the sense that   
\begin{equation}\label{f}
f\in \cP_{s,{\tail_{\rm o}}}(\d)
\end{equation}
with  tail function 
\begin{equation}\label{sidone}
\ks{\d;\g}:= \frac{4}{s}\  \log\Big(e+
 \frac{2^9}{s^{n} \g\d}\Big)\,.
\end{equation}
Let
 $\K\ge  3 \KO\ge 6$
satisfying 
\beq{porticato2}
\K^{2\nu-3n-3}
\geq e^{s+5} 
 2^{n+11} n^{2n}   \frac{L}{s^{2n+1}}\, \frac1{\g\d}
 \qquad
 \mbox{for some}\quad
 \nu\ge \frac32 n+2
 \eeq
and where $L$ defined in \eqref{Lip}.
Set
\begin{equation}\label{islandaTE}
r_k:=\sqrt\e\frac{\K^\nu}{L|k|}\,.
\end{equation}
Finally assume that
\beq{dimenticavo}
\e\le \frac{(Lr)^2}{\K^{2\nu}}\,.
\eeq
Then,
for any   $k\in \genKO$ with
$\ks{\d;\g}\le \noruno{k}\le \KO$, 
there exists $\sa^{(k)}\in[0,2\pi)$ and  a symplectic change of variables defined in a neighbourhood of the simple resonance $\Dk\times \torus^n$ such that the following holds:
\begin{equation}\label{canarinoTE}
\Psi_k: 
\Dk_{r_k/2}
\times \T^n_{s(1-1/\K)^2} \to 
\Dk_{r_k} \times \T^n_{s(1-1/\K)} 
\,,
\end{equation}
and  
\begin{equation}\label{HkTE}
\boxed{
\ham\circ\Psi_k
=:
h(y)+
2|f_k|\e
\Big( \cos(k\cdot x +\sa^{(k)})+
\mathtt G^k(y,k\cdot x)+
\mathtt f^k (y,x)
 \Big)
 }
\end{equation}
where  $\mathtt G^k(y,\cdot )\in\hol_2^1$
for every $y\in \Dk_{r_k/2}$ and 
\beq{martinaTE}
\norma \mathtt G^k\norma_{\Dk, r_k/2,2}\le \g \ .
\eeq
Finally, 
\begin{equation}\label{lothlorienTE}
\proiezione_{k\Z}\mathtt f^{k}=0\qquad 
\mbox{and}\qquad 
\norma \mathtt f^{k} \norma _{\Dk,r_k/2,s(1-1/\K)/2} 
\leq
\frac{2^{10 n} n^{3n}}{ s^{3n}\d} 
e^{-\K s/8}
\ .
\end{equation}
\end{theorem}

\rem
\label{kuno}
(i)
Recalling \eqref{canuto}, estimate
\eqref{martinaTE}
means
\beq{martinaTEbis}
\sup_{y\in \Dk_{r_k/2}}
\sum_{j\in\Z}
| \mathtt G^k_j(y) | e^{2|j|} \le \g \ .
\eeq
This implies that
for every $y\in \Dk$
the $2\pi$-periodic real function 
\beq{dolore}
\sa\mapsto \cos(\sa+\sa^{(k)})+\mathtt G^k(y,\sa)
\eeq
behaves like a cosine in the sense that {\sl it 
is a Morse function with 
only one  maximum and 
one  minimum
and no other critical points}.
To prove this,  notice that by \eqref{martinaTEbis}
we have 
$$
\sup_{y\in\Dk, x\in\torus^1}
|\partial_\sa\mathtt G^k(y,\sa)|\leq \g/e^2\,,
\qquad
\sup_{y\in\Dk, x\in\torus^1}
|\partial^2_{\sa \sa}\mathtt G^k(y,\sa)|\leq \g/e^2
\,.
$$
Therefore, denoting
by $\psi(\sa)$ the derivative of the function
in \eqref{dolore}, we have that
$\psi(\sa)>0$ for 
$\sa\in(-\sa^{(k)}+\sa_*,\pi-\sa^{(k)}-\sa_*)$ 
and  $\psi(\sa)<0$ for 
$\sa\in(\pi-\sa^{(k)}+\sa_*,2\pi-\sa^{(k)}-\sa_*),$
where $\sa_*:=\arcsin(\g/e^2)$.
Moreover in the interval
$(-\sa^{(k)}-\sa_*,-\sa^{(k)}+\sa_*)$ the function
 $\psi(\sa)$ has a zero  and is strictly increasing
 since $\psi'(\sa)\geq \sqrt{1-\g/e^2}-\g/e^2=:c>0.$
 Finally 
 in the interval
$(\pi-\sa^{(k)}-\sa_*,\pi-\sa^{(k)}+\sa_*)$ it has a zero  and is strictly decreasing
 since $\psi'(\sa)\leq -c<0.$

\nl
(ii) As a consequence the phase portrait of the effective Hamiltonian
\beq{cipsters}
h(y)+
2|f_k|\e
\big( \cos(k\cdot x +\sa^{(k)})+
\mathtt G^k(y,k\cdot x) \big)
\eeq
and that of the Hamiltonian
$
h(y)+
2|f_k|\e
 \cos(k\cdot x +\sa^{(k)})
$
are topologically equivalent.
\erem

\nl
(iii) As well know the effective Hamiltonian \equ{cipsters} 
is an integrable system as it depends only on one angle.
Indeed, fix $k\in \integer^n\backslash\{0\}$ with gcd$(k_1,\ldots,k_n)=1$, 
 then, there exists a matrix 
$A_k\in\ {\rm Mat}_{n\times n}(\Z)$ such that\footnote{Here, $k$ is  a row vector.}
\begin{equation}\label{scimmia}
A_k=\binom{\hat A_k}{k}\in\ {\rm Mat}_{n\times n}(\Z)\,,\ \ \ 
\hat A_k\in{\rm Mat}_{(n-1)\times n}(\Z)\,, \ \ \
\det A_k=1\,,\ \ \ 
|\hat A_k|_{{}_\infty}\leq |k|_{{}_\infty}\ ,
\end{equation}
where $|\cdot|_\infty$ denotes the sup--norm
of the matrix and of the vector, respectively.
The existence of such a matrix is guaranteed by an elementary result of linear algebra based on Bezout's Lemma (see Lemma~\ref{alfonso} in Appendix~\ref{appendicedialfonso}).
\\
Let us perform the linear symplectic change of variables
\begin{equation}\label{talktothewind}
\Phi_k: (Y,X)  \mapsto  (y,x):=
(A_k^TY,  A_k^{-1} X)\,,
\end{equation}
which is generated by the generating function $S(Y,x):=Y\cdot A_k x$.
Note that $\Phi_k$ does not mix actions with angles, its projection on the angles is a diffeomorphism of $\T^n$ onto $\T^n$,
and, most relevantly, $X_n=k\cdot x$
is the ``secular angle''.
\\
In the $(Y,X)$--variables, the secular Hamiltonian in \equ{cipsters} takes the form
\beq{finlandia}
\mathtt h(Y)+
2|f_k|\e
\Big( \cos(X_n +\sa^{(k)})+
\mathtt G^k(A_k^TY,X_n)
 \Big)\,,\qquad
 {\rm with}\qquad
 \mathtt h(Y):=h(A_k^TY)\,.
\eeq
Fix $y_0\in\Dk$ on the exact resonance, namely $\partial_y h(y_0)\cdot k=0.$ 
Let $Y_0$ be such that
$y_0=A_k^TY_0.$
We have
$$
\partial_{Y_n} \mathtt h(Y_0)
\stackrel{\equ{scimmia}}=
\partial_y h(A_k^TY_0)\cdot k=
\partial_y h(y_0)\cdot k=0\,,\qquad
\partial^2_{Y_n Y_n} \mathtt h(Y_0)
=
\partial^2_{yy} h(y_0)k\cdot k\,,
$$
where $\partial^2_{yy} h$ is the Hessian matrix of $h$.
By Taylor expansion the secular Hamiltonian in
\equ{finlandia}
takes the form (up to an addictive constant)
\beq{melette}
\frac12 \Big(\partial^2_{yy} h(y_0)k\cdot k\Big)
(Y_n-Y_{0n})^2+
O\Big( 
(Y_n-Y_{0n})^3
\Big)
+
2|f_k|\e
\Big( \cos(X_n +\sa^{(k)})+
\mathtt G^k(A_k^TY,X_n)
 \Big)\,.
\eeq

\nl
(iv) In particular if the Hamiltonian $h$ is {\sl convex}
the coefficient $\partial^2_{yy} h(y_0)k\cdot k=:m_k$
is bounded away from zero
and the phase portrait of the secular  Hamiltonian in \equ{melette}
is topologically equivalent, for $|Y_n-Y_{0n}|$ small\footnote{Namely
in the region $\Dk\times\torus^n$ in the original variables.}, to that of the pendulum
$$
\frac12 m_k (Y_n-Y_{0n})^2 + 2|f_k|\e
\cos(X_n +\sa^{(k)})\,.
$$

\Giu

\subsubsection*{The $y$-dependent case.}

Let us briefly turn to the $y$-dependent case.
First we note that it can happen that, even if the potential $f(y,x)$ satisfies the  non-degeneracy condition
given in Definition \ref{orcadi} 
at some point $y_0$,   there is no neighborhoud of $y_0$
on which the non-degeneracy condition holds.
For example consider the potential
$$
f(y,x)=f(y_1,x):=
\frac12 \sum_{k\neq 0} \left(\noruno{k}^{-n}-\frac{y_1}{r}\right)
e^{-\noruno{k}s} e^{\ii k\cdot x}\,.
$$
We have that $\|f\|_{D,r,s}=1$ with $D=\{0\}$
and
$$
f(0,\cdot)\in  \cP_{s,0}(1/2)\,.
$$
However, $f_k(r j^{-n})=0$ for every
$|k|=j;$
in particular for every odd number $j=2h+1,$ $h\geq 1$ and
$k:=(h+1,h,0,\ldots,0)\in\gen.$
Then for every $\d>0$, tail function $\tau>0$ and odd $j\geq 3,$
we have
$$
f(r j^{-n},\cdot)
\notin \cP_{s,\tau}(\d)\,.
$$ 
Then,
we will prove that the non-degeneracy condition holds
in a set of large measure.
In particular
 we fix $\mu>0$ and prove that, if for a certain point
$y_0\in D$ the potential $f(y_0,\cdot)\in \cP_{s,\tau_*}(\d)$
for a suitable $\tau_*=\tau_*(\mu)$ and 
$|k|\geq \tau_*(\mu),$
 then
\eqref{HkTE}-\eqref{lothlorienTE} holds
(with $f_k=f_k(y)$ and for a suitable phase $\sa^{(k)}=\sa^{(k)}(y)$)
 {\sl for every $y\in B_{r/2e}(y_0)$  up to a set of relative measure
 smaller that $\mu$}.
 \\
For the  precise statement, we refer to   Theorem~\ref{scoziaY} below.

\section{Functional setting and generic holomorphic classes} 

%

\subsection{Analytic function spaces}
\label{spigola}

\begin{itemize} 

\item[\bf (a)]
{\bf Norms and neighbourhoods} 

\nl
In this paper $|\cdot|$
 denotes the standard Euclidean norm on $\complex^n$ and its subspaces, namely 
 $|x|^2=\sum_{j=1}^n|x_j|^2$ for every $x\in\complex^n$. 

\nl
 $\noruno{k}$  denotes the 1-norm $\sum |k_j|$.

\nl
For linear maps and matrices $A$ (which we shall always identify),
$|A|$ denotes the  standard ``operator norm'' $\dst |A|=\sup_{u\neq 0} |Au|/|u|$. 

\nl
$|M|_{{}_\infty}$, with $M$ matrix (or vector), denotes the maximum norm $\max_{ij}|M_{ij}|$ (or $\max_i |M_i|$).

\nl
We shall use three different (non--equivalent) norms on holomorphic functions with domain $ \T^n_s$, $D_{r}\times \T^n_s$
or $D_{r}$ ($D$ being a susbset of $\R^n$):
given a  holomorphic function $f$ with values in $\C^m$ and  domain  $ \T^n_s$, $D_{r}\times \T^n_s$ or $D_{r}$, with $s,r> 0$ we denote by 
$$\dst \sum_{k\in\Z^n} f_k e^{\ii k\cdot x}\ \quad {\rm or}\quad \dst f(y,x)=\sum_{k\in\Z^n}f_k(y) e^{\ii k\cdot x}$$
its Fourier expansion and define the following {\sl sup--norm, $\ell^\io$--Fourier norm and $\ell^1$--Fourier norm}:

\beqa{ionicoionico}
&&{\modulo}f{\modulo}_s:=\sup_{\T^n_s}|f|\ ,
\quad
{\modulo}f{\modulo}_{r}:=\sup_{D_r}|f|\ ,
\quad
{\modulo}f{\modulo}_{r,s}:=\sup_{D_r\times \T^n_s}|f|\ ,
\\
&&\label{enorme}
\|f\|_s:=\sup_{k\in \integer^n} |f_k| e^{|k|_{{}_1}s}\ ,\quad 
\|f\|_{r,s}
 := \sup_{k\in\Z^n}
 \big(
 \sup_{y\in D_r}|f_k(y)|e^{|k|_{{}_1}s}
 \big)
\ ,\\
&&\label{attila}
\norma f\norma_s:=
\sum_{k\in\Z^n} 
|f_k| e^{|k|_{{}_1}s}\ ,\quad
\norma f\norma_{r,s}:=
\sup_{y\in D_r}\sum_{k\in\Z^n} 
|f_k(y)| e^{|k|_{{}_1}s}\ .
\eeqa
If the (real) domain need to be specified, we set, respectively, 
\beq{zenzero}
|f|_{D,r,s}:=|f|_{r,s}\ ,\quad \|f\|_{D,r,s}:=\|f\|_{r,s}\ ,\quad  \norma f\norma_{D,r,s}:=\norma f\norma_{r,s}\ .
\eeq
\end{itemize}

\rem\label{ozono}
(i) The space of functions $f:\T^n_s\to \complex^m$ endowed with the sup--norm $|\cdot|_s$ 
or the  $\ell^1$--Fourier norm $\norma \cdot \norma_{s}$
is a Banach algebra, while $\{f:\T^n_s\to \complex^m \st \|f\|_s<\io\}$ is just a Banach space (not a Banach algebra).  However, the norm $\|\cdot\|_s$ is particularly suited to describe $\{f:\T^n_s\to \complex\st \|f\|_s<\io\}$ as a probability space; compare item~{\bf (h)} below.

\nl
(ii)  As already mentioned the three norms in \equ{zenzero}  are not equivalent.  Indeed,  for
any $\s>0$, one has\footnote{We have
$\dst \sum_{k\in\Z^n\setminus 0}e^{-|k|_{{}_1}\s}
= \coth^n(\s/2)-1.$
Moreover $\coth^n x -1\leq (n/x)^n.$ Indeed for $0<x\leq 1$
the estimates follows by 
$\coth x<2/\sinh x<2/x.$
In the case $x>1$ we have 
$$
\coth^n x -1\leq   (1+e^{1-2x})^n-1\leq n(1+1/e)^{n-1}e^{1-2x}\leq  (n/x)^n\,,
$$
where in the second inequality we have used that 
$(1+y)^n\leq 1+n (1+1/e)^{n-1} y$ for $0\leq y\leq 1/e$, while in the last one
we exploit $\max_{x\geq 1}x^n e^{-2x}=(n/2e)^n.$
}
\beqa{battiato3}
\|f\|_{r,s}\leq {\modulo}f{\modulo}_{r,s}\leq 
\norma f\norma_{r,s}&\le& 
(\coth^n(\s/2)-1)\|f\|_{r,s+\s}\\
&\leq&
(2n/\s)^n\|f\|_{r,s+\s}\,.\nonumber
\eeqa
\erem

\giu
The Banach subspace of $\{f:\T^n_s\to \complex:\|f\|_s<\io\}$ of real--analytic functions with zero average ($f_0=0$) will be the natural ambient function space. Generic  elements of  such a space will be the typical potentials   to which our uniformaveraging theory applies. We givet it a name:

\nl
The following two definitions  are needed in order to decompose a holomorphic function on $\torus^n$ into a sum over generators of 1-d maximal lattices of  holomorphic functions on $\torus^1$. Later the Fourier modes $k\in\gen $ will be identified with simple resonances.

\begin{itemize}
\item[\bf (b)] {\bf Lattice Fourier projectors}

\nl
Given $f(y,x)=\sum_{k\in\Z^n}f_k(y)e^{\ii k\cdot x}$ and   a sublattice $\L$ of $\Z^n$, 
we denote by $\proiezione_\L$ the projection on the Fourier coefficients in $\L,$ namely
\beq{normandia}
\proiezione_\L f:=\sum_{k\in\L}f_k(y)e^{\ii k\cdot x}\,.
\eeq
and by $\proiezione_\L^\perp$ its ``orthogonal'' operator (projection on the Fourier modes in $\integer^n\bks\L$):
$$
\proiezione_\L^\perp f:=\sum_{k\notin\L}f_k(y)e^{\ii k\cdot x}\,.
$$
Obviously
\begin{equation}\label{galles}
\norma \proiezione_\L f\norma _{r,s}\,, \ \norma \proiezione_\L^\perp f\norma _{r,s}\leq \norma f\norma _{r,s}\,.
\end{equation}

\item[\bf (c)] {\bf  Fourier Truncation operators}

\nl
Given $N>0$, we introduce   the following ``truncation'' and ``high--mode'' operators $T_N$ and $\TNp$:
\begin{equation}\label{labestia}
T_N f (y,x):=\sum_{|k|_{{}_1}\leq N}f_k(y)e^{\ii k\cdot x}\,,
\qquad
\TNp f (y, x):=\sum_{|k|_{{}_1}>N}f_k(y)e^{\ii k\cdot x}
\,.
\end{equation}
Note that $\proiezione_\L$ and $T_N$ commute.
Note, also,  that
\begin{equation}\label{galles2}
\norma T_N f\norma _{r,s}\,, \ \norma \TNp f\norma _{r,s}\leq \norma f\norma _{r,s}\,,
\end{equation}
and  that 
\begin{equation}\label{582}
\norma \TNp f\norma _{r,s-\s}\leq e^{-(N+1)\s} \norma f\norma _{r,s}\,,\qquad 0<\s<s\,.
\end{equation}

\end{itemize}

%

\begin{itemize}

\item[\bf (d)] 
{\bf 1d--Fourier  projectors} 
 
\nl 
 Given $f\in  \hol_s^n$ and $k\in \gen$  we define the following {\bf 1d--Fourier projector}
\beq{pale}
f\in \hol_s^n \mapsto \pi_{k\integer} f =:F^k\in\hol_{|k|_{{}_1}s}^1 \quad {\rm where}\quad
F^k(\sa):=\sum_{j\in \integer\backslash\{0\}} f_{jk}e^{\ii j \sa}\ ,
\eeq
$f_{jk}$ being the Fourier coefficient of $f$ with Fourier index $jk\in\integer^n$.
\end{itemize}

\nl
It is immediate to see that:

\nl
{\sl  Any $f\in \hol_s^n$ can be  uniquely written as:}
\beq{dec}
\boxed{f(x)= \sum_{k\in \gen} F^k(k\cdot x)}
\eeq 
Notice also that,  if $k\in \gen$ and $\norma f \norma_{r,s}<\io$, then
\begin{equation}\label{stanco}
\norma F^k\norma_{r, |k|_{{}_1}s}
\leq 
\norma f \norma_{r,s}\,.
\end{equation}


%

\noi
We are now ready to define the main  function spaces.

\giu
Finally, we introduce  a  probability measure on the unit ball in $\hol^n_s$.

\begin{itemize}
\item[\bf (e)]
Denote by
$\ells$  the Banach space of complex sequences (over $\Z^n_\varstar$)  given by 
\beq{diciamolo}
\ells:=\big \{z\in \complex^{\Z^n_\varstar}
\st z_k\neq 0 \ {\rm and}\ 
|z|_{{}_\infty}:=\sup_{k\in \integer^n_\varstar } |z_k|<+\io\big\}\ .
\eeq
Then, the map 
\begin{equation}\label{miserere}
j:f\in  \hol_s^n \to \big\{f_k e^{|k|_{{}_1}s}\big\}_{k\in \integer^n_\varstar}\in \ells
\end{equation}
is an isomorphism of Banach spaces\footnote{Recall that since the functions in $\hol^n_s$ are {\sl real}--analytic one has the reality condition $f_k=\bar f_{-k}$.},
which allows to identify functions in $ \hol_s^n$ with points in $\ells$ and the Borellians  of $ \hol_s^n$ with those of $\ells$. 
\\
Denote by $\palla$ the closed ball of radius one in $\hol_s^n$ and by $\cB$ the Borellians in $\palla$.
\\
On $\palla$ we can introduce  the following natural (product) {\sl probability measure}.\\ 
Consider, first,  the probability measure given by the  normalized Le\-bes\-gue--product measure 
on the unit closed ball of $\ells$, namely,   the unique probability measure $\mu$ on the Borellians of $\{z\in \ells: |z|_{{}_\infty}\le 1\}$ such that,
given   Lebesgue measurable  sets $\Ak $ in the unit complex disk $\Ak \subseteq D:=\{w\in\complex:\ |w|\le 1\}$
with $\Ak \neq D$ only for finitely many $k$, one has
 $$
\mu \Big(\prod_{k\in \integer^n_\varstar} \Ak \Big)=
\prod_{\{k\in \integer^n_\varstar:\, \Ak \neq D\}} 
\frac{1}{\pi}{\rm meas}(\Ak )\,
$$
where ``meas'' denotes  the Lebesgue measure on the unit complex disk $D$. 
\\
Then,  {\sl the isometry $j$ in \eqref{miserere}
naturally induces a} {\bf  probability measure\footnote{I.e, $\mu_s(\palla)=1$.}
$\mu_s$ on the Borellians $\cB$}.

\end{itemize}

\subsection{Generic properties  of periodic holomorphic classes} 
\label{bretagna}

Here we discuss some properties of the classes $\cP_{s,\tail}$
of non-degenerate 
introduced in Definition \ref{orcadi}.

\rem\label{pioggia} (i) 
Since $f\in \hol^n_s$, one has that $|f_k|\le \|f\|_s e^{-\noruno{k} s}$ for all $k$'s and \equ{star} says that, when $k$ is a generator of maximal 1d--lattices (later corresponding to simple resonances), the $k$--Fourier coefficient does not vanish and  is controlled in a quantitive way from below: $|k|_{{}_1}^{-n}$ is a suitable weight (needed in the proof of Proposition~\ref{grizman} below), while $\d$ is any number satisfiying
\beq{bha}
 \inf_{\noruno{k}>\tail(\d)} |f_k|  |k|_{{}_1}^{{n}}\ e^{|k|_{{}_1}s} \ge \d>0\ .
\eeq

\nl
(ii) 
It is easy to construct functions in $\cP_{s,\tail}(\d)$. For 
example let
\beq{esempietto}
 f(x):=2\d \sum_{k\in \gen} 
|k|_{{}_1}^{-{n}} e^{-|k|_{{}_1}s}\,  \cos (k\cdot x)\ ,
\eeq
which has Fourier coefficients
\beqno
 f_k=\casialt{\dst\d |k|_{{}_1}^{-{n}} e^{-|k|_{{}_1}s}, }{\pm k\in\gen}{0 ,}
\eeqno
and 1d--Fourier projections 
$$
 F^k(\sa)=\d |k|_{{}_1}^{-{n}} e^{-|k|_{{}_1}s} \cos \sa\ .
$$
Then,  $ f\in \cP_{s,0}(\d)$ and, also,    $ f\in \cP_{s,\tail}(\d)$ for any  choice of tail function $\tail(\delta)$.
\erem

\noi
Here, we show that the classes $\cP_{s,\tail}$ are ``general'' in several (topological and measure theoretical) ways.

\begin{proposition}\label{grizman}{\bf (Properties of $\cP_{s,\tail}$)}
Let $s>0$ and $\tail$ be a tail function. Then:
\begin{itemize}
\item[{\rm (i)}] 
The set $\cP_{s,\tail}\subseteq \hol_s^n$ contains an open dense set.

\item[{\rm (ii)}]
$\cP_{s,\tail}\, \cap \, \palla\in \cB$ and $\m_s(\cP_{s,\tail}\cap \palla)=1$.

\item[{\rm (iii)}]
$\cP_{s,\tail}$  is a  prevalent set\footnote{We recall that a Borel set $P$ of a Banach space $X$ is called {\sl prevalent} if there exists a compactly supported probability measure $\nu$ on the Borellians of $X$  such that $\nu(x+P)=1$ for all $x\in X$; compare, e.g., 
\cite{HK}
\label{nurzia}}.

\end{itemize}
\end{proposition}

\proof

\giu
(i)  {\sl  $\cP_{s,\tail}$ contains an open subset $\cP_{s,\tail}'$ which is dense in the unit ball of  $ \hol_s^n$.}

\nl
Let us define $\cP_{s,\tail}'$  as $\cP_{s,\tail}$ but with the difference that  \equ{star} is replaced by the
{\sl stronger condition}\footnote{Note, however,  that $\m_s(\cP'_s)=0$.}

\beq{starstar}
\exists\, \d>0\ {\rm s.t.} 
\ |f_k|\geq\d\ e^{-|k|_{{}_1}s} \,,\ \ 
\forall\ \,k\in\gen,\  |k|_{{}_1}> \tail(\delta)
\eeq

\nl
Let us first prove that $\cP_{s,\tail}'$ is open. Let $f\in\cP_{s,\tail}'$.
We have to show that there exists $\rho>0$ such that if $\|g\|_s< \rho$, then $f+g\in\cP_{s,\tail}'$.
Fix $\d>0$ such that \equ{starstar} holds and, by continuity of $\tail(\delta)$, choose $\rho<\d$ small enough such that
$
[\tail(\delta)] > \tail(\delta')-1\,$, where $ \d':=\d-\rho
$
and $[\cdot]$ denotes integer part.
Then, since $\tail(\delta)$ is not increasing,
 it is immediate to verify that $|k|_{{}_1}> \tail(\delta) \iff |k|_{{}_1}> \tail(\delta')$.
 Moreover
$$
|f_k+g_k|e^{|k|_{{}_1}s}\geq |f_k| e^{|k|_{{}_1}s} -\|g\|_s \geq \d-\rho=\d'\,, \qquad
\forall\ \,k\in\gen,\  |k|_{{}_1}> \tail(\delta')\,,
$$
namely $f+g$ satisfies  \equ{starstar} (with $\d'$ instead of $\d$).

\nl
Let us now show that $\cP_{s,\tail}'$  is dense in the unit ball of  $ \hol_s^n$.
Take $f$ in the unit ball of  $ \hol_s^n$ and $0<\loge <1$. We have to find $\tilde f\in\cP_{s,\tail}'$
with $\|\tilde f-f\|_s\leq \loge $.
Let $\d:=\loge /4$ and denote by  $f_k$ and $\tilde f_k$ (to be defined) be the Fourier coefficients of, respectively,  
$f$ and $\tilde f$. We, then, let $\tilde f_k=f_k$ unless
 $k\in\gen$, $|k|_{{}_1}> \tail(\delta)$ and $|f_k|e^{|k|_{{}_1}s}< \d$, in which case, $\tilde f_k= \delta e^{-|k|_{{}_1}s}$.
It is, now,  easy to check that  $\tilde f\in \cP'_s$ and is $\loge $--close to $f$. 

\Giu
(ii) {\sl $\cP_{s,\tail}\, \cap \, \palla\in \cB$ and $\m_s(\cP_{s,\tail}\cap \palla)=1$}

\giu
We shall prove that, for every $\d>0$, 
the measure of the sets of
potentials $f$ that do not satisfy  \equ{star} is
$O(\d^2)$, the result will follow letting $\d\to 0$.

\nl
By the identification \eqref{miserere}, the
measure of the set of
potentials $f$ that do not satisfy \equ{star}
with a given $\delta$
 is  bounded by
 \beq{babbonatale}
  \d^2\, \sum_{k\in\integer^n} |k|_{{}_1}^{-2{n}}\,.
  \eeq
\begin{remark}\label{befana}
Recalling footnote \ref{befana2}, one could impose the condition
$|f_k|\geq \d |k|_{{}_1}^{-\npiu}\ e^{-|k|_{{}_1}s}$
in \equ{star}. Then  \eqref{babbonatale}
 would become
 $ \d^2\, \sum_{k\in\integer^n} |k|_{{}_1}^{-2\npiu},$
 which is still fine if $\npiu>n/2.$
\end{remark}

\Giu
(iii) {\sl  $\cP_{s,\tail}$ is prevalent}.

\giu
Consider the following compact subset of
$\ells$: 
let $\mathcal K:=\{ z=\{z_k\}_{k\in \integer^n_\varstar  } : z_k\in D_{1/|k|_{{}_1}} \}$,
where $D_{1/|k|_{{}_1}}:=\{w\in\complex:\ |w|\le 1/|k|_{{}_1}\}$, 
and let  $\nu$ be 
the unique probability measure supported on $\mathcal K$ such that,
given  Lebesgue measurable sets 
$\Ak \subseteq D_{1/|k|_{{}_1}}$,
with $\Ak \neq D_{1/|k|_{{}_1}}$
only for finitely many $k$, one has
$$
\nu \Big(\prod_{k\in \integer^n_\varstar} \Ak \Big):=
\prod_{\{k\in \integer^n_\varstar:\, \Ak \neq D_{1/|k|_{{}_1}}\}} 
\frac{|k|_{{}_1}^2}{\pi}{\rm meas}(\Ak )\,.
$$ 
The isometry $j_s$ in \eqref{miserere}
naturally induces a  probability measure
$\nu_s$ on  $\hol^n_s$ with support in the compact set $\mathcal K_s :=j_s^{-1}\mathcal K$.
Reasoning as in the proof of 
$\mu_s(\cP_{s,\tail})=1$, 
one can show that
$\nu_s(\cP_{s,\d})
\geq 1-{\rm const}\, \d^2$.
It is also easy to check that, for every $g\in  \hol_s^n$, the translated 
set $\cP_{s,\d}+g$ satisfies
$\nu_s(\cP_{s,\d}+g)\geq \nu_s(\cP_{s,\d})$.
Thus,  one gets 
$\nu_s(\cP_{s,\tail}+g)= \nu_s(\cP_{s,\tail})
= 1$, $\forall\ \, g\in  \hol_s^n$, which means  that
$\cP_{s,\tail}$ is prevalent
(recall footnote~\ref{nurzia}). 
\qed



\section{A normal form lemma with ``small'' analyticity loss}\label{vongole}
In this section we describe an analytic normal form lemma for nearly--integrable Hamiltonians $H(y,x)=h(y)+f(y,x)$,
which  allows to average out non--resonant Fourier modes of the perturbation $f$ on suitable non--resonant regions,
and allows for ``very small'' analyticity loss in the angle variables, a fact, which will be crucial in our applications.

\nl
We recall (\cite{nek}, \cite{poschel}) that, given an integrable Hamiltonian $h(y)$, positive numbers $\a,K$ and a lattice $\L\subset \Z^n$, a (real or complex) domain $U$ is $(\a,K)$ non--resonant modulo $\L$ (with respect to $h$) if  
\begin{equation}\label{panzapiena}
|h'(y)\cdot k |\ge \a\ ,
\ \ \ \forall\ y\in U\ , \forall\  k\in \integer^n\bks \L\ , \    \noruno{k}\leq K\ .
\end{equation}

\nl
The main point of the following ``Normal Form Lemma'' is that the ``new'' averaged Hamiltonian is defined, 
in the fast variable (angle) domain, in a region ``almost equal'' to the original domain, ``almost equal'' meaning a
complex strip of width $s(1-1/K)$ if $s$ is the width of the initial  angle analyticity. More precisely, we have:

\begin{proposition}[Normal form with ``small'' analyticity loss]\label{pesce}\ \\
Let $r,s,\a>0$,  $K\in\N,$ $K\ge 2$, $D\subseteq \real^n$,  and let $\L$ be a lattice of $\integer^n$.
Let 
\begin{equation}\label{pegaso}
H(y,x)=h(y)+f(y,x)
\end{equation}
 be real--analytic on  
$ D_r \times \T^n_s$ with $\norma f\norma _{r,s}<\infty.$
Assume that $D_r$ is
($\a$,$K$) non--resonant modulo
$\L$ 
 and that
\begin{equation}\label{gricia}
\tetta_\varstar := \frac{2^{11} K^2}{\a r s}\, \norma f\norma _{r,s}  <1\,.
\end{equation}
Then,  there exists a real--analytic 
 symplectic change of variables
\begin{equation}\label{dadostar}
\Psi: (y',x')\in D_{r_\varstar}\times \T^n_{s_\varstar}\  \mapsto\ 
(y,x)\in D_r \times \T^n_s \, \quad {\rm with} 
\quad
r_\varstar:=r/2\,,\ \ 
s_\varstar:=s(1-1/K)
\end{equation}
satisfying
\begin{equation}\label{argo0}
\noruno{y-y'}\leq \frac{\tetta_\varstar}{2^7 K} r\,,\qquad
\max_{1\leq i\leq n}|x_i-x'_i|
\leq 
\frac{\tetta_\varstar}{16 K^2} s\,,
\end{equation}
and such that 
\beq{senzanome}
H\circ\Psi=h+ f^\flat +
f_\varstar\,,   \qquad f^\flat := \proiezione_\L f+ \Tp\proiezione_\L^\perp f
\eeq
with 
\begin{equation}\label{pirati}
\norma f_\varstar\norma _{r_\varstar,s_\varstar}
\leq 
\frac{1}{K}\tetta_\varstar\norma f\norma_{r,s}
\,,\qquad 
\norma T_{K}\proiezione_\L^\perp f_\varstar\norma _{r_\varstar,s_\varstar} 
\leq
(\tetta_\varstar/8)^{{K}}\frac{8}{eK}\norma f\norma _{r,s}\,.
\end{equation}
Moreover,  re-writing \equ{senzanome} as 
\beq{silvia}
H\circ \Psi=h+g+ f_{\varstar\varstar}\qquad \mbox{where}\quad \proiezione_\L g=g\ ,\quad \proiezione_\L f_{\varstar\varstar}=0\,,
\eeq
one has
\beq{giovanni}
\norma g- \proiezione_\L f\norma _{r_\varstar,s_\varstar}\le 
\frac{1}{K}\tetta_\varstar \norma f\norma _{r,s}\ ,\qquad
\norma f_{\varstar\varstar}\norma _{r_\varstar,s/2}
\le 
 2 e^{-(K-2)\bs}\norma f\norma _{r,s}\,,
\eeq
where
\begin{equation}\label{shreck}
\bs:=\min\left\{
\frac{s}{2},\,
\ln\frac{8}{\tetta_\varstar}
\right\}\,.
\end{equation}

\end{proposition}

\rem\label{cb500x} $\phantom.$
\begin{itemize}
\item[{\rm (i)}] 
The ``novelty'' of this lemma is that the bounds in \equ{pirati} and the first one in \eqref{giovanni}
 hold on the large angle domain $\torus^n_{s_\varstar}$ with $s_\varstar=s(1-1/K)$. In particular
 the first estimate in \equ{pirati}
(or, equivalently, in \eqref{giovanni})
   will be important in our analysis  in order to obtain \eqref{552bis}, \eqref{cristina} 
and, therefore, \eqref{120}, \eqref{martha} and finally \eqref{martina},
which is the key to prove \eqref{martinaTE} in  Theorem 
\ref{scozia}. The drawback of the gain in angle--analyticity strip is that the power of $K$ in the smallness condition 
\equ{gricia} is not optimal: for example in \cite{poschel} the power of $K$ is one (but $s_\varstar=s/6$, which would not work in our applications). 

\item[{\rm (ii)}]
Having information on non--resonant Fourier modes up to order $K$, the best one can do is to average out
the non--resonant Fourier modes up to order $K$, namely, to ``kill" the term $T_K \proiezione_\L^\perp f$ of the Fourier expansion of the perturbation. This explains the ``flat'' term $f^\flat=\proiezione_\L f+ \Tp\proiezione_\L^\perp  f$ surviving in \equ{senzanome} and  which cannot be removed in general.
Now, think of the remainder term $f_\varstar$ as 
$$f_\varstar=\proiezione_\L f_\varstar + \big( \Tp\proiezione_\L^\perp  f_\varstar + T_{K}\proiezione_\L^\perp  f_\varstar\big)\ ;$$
then,  $\proiezione_\L f_\varstar$  is a $(\tetta_\varstar \norma f\norma _{r,s}/K)$--perturbation of the part in normal form (i.e., with Fourier modes in $\L$), while 
$ \Tp\proiezione_\L^\perp  f_\varstar$ is, by \equ{582}, a term exponentially small with $K$ (see also below)
and 
$T_{K}\proiezione_\L^\perp  f_\varstar$ is a very small remainder bounded by $8(\tetta_\varstar /8)^K \norma f\norma _{r,s}/eK$.

\item[{\rm (iii)}] 

We note that 
\eqref{silvia} follows from \equ{senzanome}. 
Indeed  we take
$$
g=\proiezione_\L f + \proiezione_\L f_\varstar\,,\qquad
f_{\varstar\varstar}=\Tp\proiezione_\L^\perp f+\proiezione_\L^\perp f_\varstar= 
T_K\proiezione_\L^\perp f_\varstar+ \Tp\proiezione_\L^\perp (f_\varstar+f)\,.
$$ 
Then the first estimate in \eqref{giovanni}
follows by the first bound in \eqref{pirati} and \eqref{galles}.
Regarding the  
second estimate in
 \equ{giovanni}, 
 we first note by \eqref{pirati} and   \equ{582} (used with
 $f\rightsquigarrow f_\varstar$
  $N\rightsquigarrow K$, $r\rightsquigarrow r_\varstar$, $s\rightsquigarrow s_\varstar$,  and $\s\rightsquigarrow\frac{s}2-\frac{s}K$ so that $s_\varstar-\s=s/2$ and $e^{-(K+1)\s}\le 
  e^{-(K-2)s/2} $)
 $$
  \norma \Tp f_\varstar\norma _{r_\varstar,s/2}=
 \norma \Tp f_\varstar\norma _{r_\varstar,s_\varstar-\s}
 \leq
 e^{-(K+1)\s}\norma f_\varstar\norma _{r_\varstar,s_\varstar}
 \leq 
  e^{-(K-2)s/2}
 \tetta_\varstar \norma f\norma _{r,s}/K\,.
 $$ 
 By \eqref{galles}, \eqref{pirati} and \eqref{582} we get
\beqano 
\norma f_{\varstar\varstar}\norma _{r_\varstar,s/2} 
&\le & 
\norma T_K\proiezione_\L^\perp f_\varstar\norma _{r_\varstar,s/2} 
 + \norma \Tp f_\varstar\norma _{r_\varstar,s/2}
+\norma \Tp f\norma _{r_\varstar,s/2}
\\
&\le& 
(\tetta_\varstar/8)^{{K}}\frac{8}{eK}\norma f\norma _{r,s}
+ e^{-(K-2)s/2}   
(\tetta_\varstar/K 
+ e^{-3s/2}) \norma f\norma _{r,s} 
\\
&\le  & 2 e^{-(K-2)\bs} \norma f\norma _{r,s}\ .
\eeqano

\item[{\rm (iv)}] 
Let us compare our results with more standard formulations, such as the Normal Form Lemma in \S~2 of \cite{poschel}.  
In that formulation, imposing the weaker  smallness condition
$\norma f\norma_{r,s}\leq \const \a r /K,$
the normal form Hamiltonian
writes
$h+\mathtt g +\mathtt f$ with $\mathtt f$ exponentially small
(of order $\norma f\norma_{r,s} e^{-Ks/6}$)
and, regarding $\mathtt g$ one knows that
\begin{equation}\label{pescetto}
\norma
\mathtt g -T_K\proiezione_\L f
\norma_{r/2,s/6}
\leq
{\rm const.}\frac{K}{\a r}
\norma f\norma_{r,s}^2\,.
\end{equation}
For our purposes  
we need to prove that, when $k\in\Z^n_\sharp,$
$|k|_{{}_1}\leq 	\KO\leq K$  ($k\in \Z^n_{\KO}$ indexes the simple resonance we want to consider while
$l\in \Z^n_{K}$ indexes the second order resonance beyond
$k$) and 
$|f_k|/\norma f\norma _{r,s}\geq \d |k|_{{}_1}^{-{n}}e^{-|k|_{{}_1}s},$ 
the quantity
$$
\frac{1}{|f_k|}\sup_{y\in D_{r/2}}
| g_k(y)-f_k|
$$
is small.
Indeed by \eqref{giovanni} we have
\begin{eqnarray}\label{albione}
&&\frac{1}{|f_k|}\sup_{y\in D_{r/2}}
| g_k(y)-f_k|
\leq
\frac{\tetta_\varstar}{K} \norma f\norma _{r,s} \frac{e^{-|k|_{{}_1}s_\varstar}}{|f_k|}
\leq
\frac{\tetta_\varstar}{K} \frac{|k|_{{}_1}^{{n}}e^{(s-s_\varstar)|k|_{{}_1}}}{\d}
\nonumber
\\
&& 
=\frac{\tetta_\varstar}{K} \frac{|k|_{{}_1}^{{n}}e^{s|k|_{{}_1}/K}}{\d}
\leq 
\frac{\tetta_\varstar}{K} e^s\frac{\KO^{{n}}}{\d}\,,
\end{eqnarray}
which is small when 
\begin{equation}\label{polenta}
\KO\ll \left(\frac{\a r s \d}{K\norma f\norma _{r,s}}\right)^{1/{n}}\,.
\end{equation}
Consider, for example, the function $f=\e \hat f$
with  $\e$ small and $\hat f$ defined in
\eqref{esempietto}.
We have that $\norma f\norma _{r,s}=c\d\e$, for a suitable constant $c>0$.
In this case  \eqref{polenta} writes
\begin{equation}\label{polenta2}
\KO\ll \left(\frac{\a r s}{K\e}\right)^{1/{n}}\,.
\end{equation}

On the other hand by estimate
\eqref{pescetto} one only have
$$
\frac{1}{|f_k|}\sup_{y\in D_{r/2}}
|\mathtt g_k(y)-f_k| \leq 
{\const. }
\frac{K\e}{\a r}  
|k|_{{}_1}^{{n}}
e^{|k|_{{}_1}s}
e^{-|k|_{{}_1}s/6}
\leq 
{\const. }
\frac{\e K}{\a r} 
K_0^{{n}}
e^{\frac56 K_0s}
\,,
$$ 
which is  small only for
\begin{equation}\label{polenta3}
\KO\ll \frac{6}{5s} \ln \frac{\a r}{K \e}\,,
\end{equation}
that is a considerably stronger bound than the one
in \eqref{polenta2}.
\\
Since we are considering simple resonances
indexed by $|k|_{{}_1}\leq \KO,$
   the non resonant  region will be non-resonant
   only up to order $\KO;$ therefore we have that the perturbation,
after normal form in the non-resonant region, will be of magnitude
$$
\e e^{-\KO s/6}\gg\e (K \e/\a r)^{1/5}\,,
$$
when the bound \eqref{polenta3} applies.
This estimate is very bad.
On the other hand, in our case, the weaker bound \eqref{polenta2} applies and we obtain that 
the perturbation is exponentially small.

\end{itemize}

\erem




\noindent
Given a function $\phi$ we denote by $X_\phi^t$
the hamiltonian flow at time $t$ generated by $\phi$
and by ``ad''  the linear operator  $u\mapsto {\rm ad}_\phi u:=\{u,\phi\}$ and ${\rm ad}^\ell$ its iterates:  
$$
{\rm ad}^0_\phi u:=u\,,  \qquad
{\rm ad}^\ell _\phi u:=\{ {\rm ad}^{\ell -1}_\phi u,\phi\}\,, \qquad \ell \geq 1\,,
$$
as standard, $\{\cdot, \cdot\}$ denotes Poisson bracket\footnote{Explicitly, 
$\dst \{u,v\}= \sum_{i=1}^n (u_{x_i}v_{y_i}- u_{y_i}v_{x_i})$.}.

\nl
Recall the identity (``Lie series expansion'') 
\begin{equation}\label{ellade}
u\circ X_\phi^1 =\sum_{\ell \geq 0}
 \frac{1}{\ell !} {\rm ad}^\ell_\phi u=
\sum_{\ell=0}^\io \frac{\partial_t^\ell   (u\circ X_\phi^t)}{\ell!} \Big|_{t=0} 
 \, ,
\end{equation}
valid for analytic functions  and small $\phi$.
We recall the following  technical lemma by \cite{poschel}.

\begin{lemma}[Lemma B.3 of \cite{poschel}]\label{argolide}
For $0<\rho<r,$ $0<\s<s,$  $D\subseteq \R^n$
$$
\sup_{y\in D_{r}}\sum_{1\leq i\leq n} \norma\partial_{x_i} \phi(y,\cdot)\norma_{s-\s}
\leq \frac{1}{e \s} \norma \phi\norma_{r,s}\,,
\qquad
\sup_{y\in D_{r-\rho}}\max_{1\leq i\leq n} \norma\partial_{y_i} \phi(y,\cdot)\norma_{s}
\leq \frac{1}{\rho} \norma \phi\norma_{r,s}\,,
$$
\end{lemma}
By Lemma \ref{argolide} we get (see also Lemma B.4 of \cite{poschel})
\begin{lemma}\label{messene}
For $0<\rho<\bar r:=\{r_0,r\},$ $0<\s<\bar s:=\{s_0,s\},$
\begin{equation}\label{tirinto}
\norma \{f,g\}\norma _{\bar r-\rho,\bar s-\s}
\leq
\frac 1e
\left(
\frac{1}{(r_0-\bar r+\rho)(s-\bar s+\s)} +
\frac{1}{(r-\bar r+\rho)(s_0-\bar s+\s)}
\right)
\norma f\norma _{r_0,s_0} \norma g\norma _{r,s}\,.
\end{equation}
\end{lemma}

\noi
Summing the Lie series in \eqref{ellade} (see  Lemma B5 of \cite{poschel}) we get, also,

\begin{lemma}\label{olimpiabis}
Let $0<\rho<r_0$ and 
$0<\s< s_0.$
Assume that
\begin{equation}\label{tebe}
\hat\tetta:=\frac{4 e \norma \phi\norma _{r_0,s_0}}{\rho\s}\leq 1\,.
\end{equation}
Then for every
$\rho<r'\leq r_0,$ $\s<s'\leq s_0,$
 the time-1-flow $X_\phi^1$
 of vector field $X_\phi$
define a good canonical transformation
\begin{equation}\label{corinto}
X_\phi^1: D_{r'-\rho}\times \T^n_{s'-\s}
\ \to \ 
D_{r'-\rho/2}\times \T^n_{s'-\s/2}
\end{equation}
satisfying
\begin{equation}\label{argo}
\noruno{y-y'}\leq \hat\tetta\frac{\rho}{4e}\,,\qquad
\max_{1\leq i\leq n}|x_i-x'_i|\leq \hat\tetta\frac{\s}{4}
\end{equation}
Moreover let $r>\rho, s>\s$ and set
$$
\bar r:=\min\{ r_0,r\}\,,\qquad
\bar s:=\min\{ s_0,s\}\,.
$$
Then for any $j\geq 0$
\begin{eqnarray}\label{leuttra}
\norma   
u\circ X^1_{\phi}-
\sum_{h \leq j} {\rm ad}^{h}_\phi u\norma  _{\bar r-\rho,\bar s-\s}
&\leq&
\sum_{h > j} \frac{1}{h !}
 \norma   {\rm ad}^{h}_\phi u\norma  _{\bar r-\rho,\bar s-\s}
 \nonumber
 \\
 &\leq&
 2(\hat\tetta/2)^j
  \norma  \{ u,\phi\}\norma  _{\bar r-\rho/2,\bar s-\s/2}
 \end{eqnarray}
for every function $u$ with $\norma u\norma _{r,s}<\infty.$
\\
In particular  when $r\leq r_0,s\leq s_0$
\begin{eqnarray}\label{delfi}
&&\norma   u\circ X_\phi^1 - u  \norma  _{r-\rho,s-\s}
\leq 
\sum_{h \geq 1} \frac{1}{h !}
 \norma   {\rm ad}^h_\phi u\norma  _{r-\rho,s-\s}
\leq 
2\hat\tetta \norma u\norma _{r,s}\,,
\\
&&\norma   u\circ X_\phi^1 - u -\{u,\phi\} \norma  _{r-\rho,s-\s}
\leq  
\hat\tetta^2 \norma u\norma _{r,s}\,,
\label{delfi2}
\end{eqnarray}
\end{lemma}
\proof
We first note that
by Lemma \ref{argolide} 
(applied with 
 $r_0\rightsquigarrow r$,
 $s_0\rightsquigarrow s$)
 for every $(y,x)\in D_{r_0-\rho}\times \T^n_{s_0-\s}$
we have
$$
\noruno{\partial_x \phi(y,x)}\leq 
\frac{1}{e\s}\norma\phi\norma_{r_0,s_0}=\frac{\hat\tetta\rho}{4e}
\leq\frac{\rho}{4e}\,,\qquad
\max_{1\leq i\leq n}|\partial_{y_i} \phi(y,x)|\leq 
\frac{1}{\rho}\norma\phi\norma_{r_0,s_0}=\frac{\hat\tetta\s}{4}
\leq\frac{\s}{4}\,.
$$
Then \eqref{corinto} holds.
\\
For $h\geq 1,$ set for brevity
$$
\norma\cdot\norma_i:=
\norma\cdot\norma_{\bar r-\frac{\rho}{2}-i \tilde\rho,\bar s-\frac{\s}{2}-i \tilde\s}\,,
\qquad
0\leq i\leq h\,,
\qquad
\tilde\rho:=\frac{\rho}{2 h}\,,\ \ \ 
\tilde\s:=\frac{\s}{2 h}\,.
$$
We get
\begin{eqnarray*}
&&\norma   {\rm ad}^{i}_\phi \{u,\phi\}\norma  _i
\\
&&
\stackrel{\eqref{tirinto}}\leq
\frac1e
\left(
\frac{1}{\tilde\rho(s_0-\bar s+i\tilde\s +\s/2)}
+\frac{1}{\tilde\s(r_0-\bar r+i\tilde\rho +\rho/2)}
\right)
\norma\phi\norma_{r_0,s_0}
\norma   {\rm ad}^{i-1}_\phi \{u,\phi\}\norma  _{i-1}
\\
&&
\leq
\frac{8h^2}{e\rho\s}\frac{1}{h+i}
\norma\phi\norma_{r_0,s_0}
\norma   {\rm ad}^{i-1}_\phi \{u,\phi\}\norma  _{i-1}\,,
\end{eqnarray*}
and, iterating,
\begin{eqnarray*}
&&\norma   {\rm ad}^{h}_\phi \{u,\phi\}\norma  _h
\leq
\frac{8h^2}{e\rho\s}\frac{h!}{(2h)!}
\norma\phi\norma_{r_0,s_0}
\norma   \{u,\phi\}\norma  _{r-\rho/2,s-\s/2}
\leq
h! (\hat\tetta/2)^h\norma   \{u,\phi\}\norma  _{r-\rho/2,s-\s/2}
\end{eqnarray*}
by Stirling's formula.
Then
$$
\sum_{h \geq j} \frac{1}{(h+1) !}
 \norma   {\rm ad}^{h+1}_\phi u\norma  _{\bar r-\rho,\bar s-\s}
\leq
\sum_{h \geq j} \frac{1}{h+1}
(\hat\tetta/2)^h\norma   \{u,\phi\}\norma  _{r-\rho/2,s-\s/2}
$$
proving \eqref{leuttra} in view of \eqref{tebe}.
\\
Finally \eqref{delfi} and \eqref{delfi2} follows by \eqref{leuttra} and
since $\norma  \{ u,\phi\}\norma  _{\bar r-\rho/2,\bar s-\s/2}\leq 2 e^{-1}\hat\tetta \norma u\norma_{r,s}$
by \eqref{tirinto}.
\eproof

\noindent
Given $K\ge 2$ and a lattice $\L$, recall the definition of $f^\flat$ in \equ{senzanome} and define 
$$
f^K:=f-f^\flat=T_{K}\proiezione_\L^\perp  f\ ,
$$
so that we have the decomposition (valid for any $f$):
\beq{decomposizione}
f=f^\flat+f^K\,,\qquad
f^\flat:= P_\L f+\Tp\proiezione_\L^\perp  f \,,\qquad
f^K:= T_{K}\proiezione_\L^\perp  f\,.
\eeq

\begin{lemma}\label{megara}
Let $0<\rho<r$ and 
$0<\s<s.$
Consider a real--analytic Hamiltonian 
\begin{equation}\label{olinto}
H=H(y,x)=h(y)+f(y,x)\qquad\mbox{analytic \ on \ } D_r\times \T^n_s\,.
\end{equation}
Suppose that $D_r$ is ($\a$,$K$) non--resonant modulo $\L$
for $h$ (with $K\ge 2$).
Assume that
\begin{equation}\label{gricia*}
\check\tetta := \frac{4 e}{\a \rho \s}\,  \norma f^K\norma _{r,s}\leq 1\,.
\end{equation}
Then there exists a real--analytic 
 symplectic change of coordinates
$$
\Psi:=X^1_\phi\, :\, D_{r_+}\times \T^n_{s_+} \ni (y',x')\
\to\ (y,x)\in D_r \times \T^n_s \,,\quad
r_+:=r-\rho\,,\ \ \ s_+:=s-\s\,,
$$
generated by a function $\phi=\phi^K=T_{K}\proiezione_\L^\perp  \phi$
with
\begin{equation}\label{maratona}
\norma \phi\norma _{r,s}\leq \norma f^K\norma _{r,s}/\a\,,
\end{equation}
satisfying
\begin{equation}\label{argo2}
\noruno{y-y'}\leq \check\tetta\frac{\rho}{4e}\,,\qquad
\max_{1\leq i\leq n}|x_i-x'_i|\leq \check\tetta\frac{\s}{4}\,,
\end{equation}
such that 
\begin{equation}\label{olintobis}
H\circ\Psi=h(y')+f_+(y',x')\,, \qquad f_+:=f^\flat+f_\varstar
\end{equation}
with
\begin{equation}\label{salamina}
\norma f_\varstar\norma _{r_+,s_+} 
\leq\ 
4\check\tetta \norma f\norma _{r,s}
\,.
\end{equation} 
\end{lemma}
Notice that, by \equ{decomposizione} and \equ{salamina}, one has
\beq{nuvole}
f_+^K=f_\varstar^K\ , \quad \norma f_+\norma _{r_+,s_+}
\le 
\norma f_\varstar\norma _{r_+,s_+}+\norma f\norma _{r,s}\le
(1+4\check\tetta)\norma f\norma _{r,s}\ .
\eeq
Notice also that 
\beq{bemolle}
f_+^\flat-f^\flat\eqby{olintobis} f_\varstar^\flat\quad \Longrightarrow \quad \norma f_+^\flat-f^\flat\norma _{r_+,s_+}\le
\norma f_\varstar\norma _{r_+.s_+}\stackrel{\equ{salamina}}\le 
 4 \check\tetta \norma f\norma _{r,s}\ .
\eeq
\proof 
Let us define
$$
\phi=\phi(y,x):=\sum_{|m|\leq K, m\notin\L}
\frac{f_m(y)}{\ii h'(y)\cdot m} e^{\ii m\cdot x}
\,,
$$
and note that $\phi$ solves the homological equation
\begin{equation}\label{tessalonica}
\{ h,\phi\}+f^K=0\,.
\end{equation}
Since $D_r$ is ($\a$,$K$) non--resonant modulo $\L$
the estimate \eqref{maratona} holds.
We now use Lemma \ref{olimpiabis} with parameters
$r_0\rightsquigarrow r, 
s_0 \rightsquigarrow s.
$
With these choices  it is  
$\hat\tetta= \check\tetta$, and, by \eqref{gricia*}  $\check\tetta\leq 1$.
Thus,  \eqref{tebe} holds and Lemma \ref{olimpiabis} applies.
\eqref{argo2} follows by \eqref{argo}.
We have
$$
H\circ\Psi=h+f^\flat+f_\varstar
$$
with
$$
f_\varstar 
=
 (h\circ\Psi-h-\{h,\phi\})+
 (f\circ\Psi-f) \,.
 $$
Since
$$
h\circ\Psi-h-\{h,\phi\}
=\sum_{\ell \geq 2}
 \frac{1}{\ell !} {\rm ad}^\ell_\phi h
 =
 \sum_{\ell \geq 1}
 \frac{1}{(\ell +1)!} {\rm ad}^\ell_\phi \{h,\phi\}
 \stackrel{\eqref{tessalonica}}=
- \sum_{\ell \geq 1}
 \frac{1}{(\ell +1)!} {\rm ad}^\ell_\phi 
f^K\,,
$$
we have
$$
\norma h\circ\Psi-h-\{h,\phi\}\norma _{r_+,s_+}
\leq 
\sum_{\ell \geq 1}
 \frac{1}{\ell!} \norma{\rm ad}^\ell_\phi 
f^K\norma_{r_+,s_+}
\stackrel{\equ{delfi}}\leq
2\check\tetta \norma f^K\norma _{r,s}
\le 
2\check\tetta \norma f\norma _{r,s}\ .$$
Finally, applying again Lemma~\ref{olimpiabis} with $u=f$, by \equ{delfi}, we get $\norma f\circ\Psi-f\norma _{r_+,s_+}
\le 2\check\tetta \norma f\norma _{r,s}$,
proving \eqref{salamina} and concluding the proof of Lemma~\ref{megara}. 
\eproof


\noi
As a preliminary step 
we apply Lemma \ref{megara}
to the Hamiltonian $H=h+f$ in \eqref{pegaso}
with $\rho=r/4$ and $\s=s/2K.$
By \eqref{galles}, \eqref{galles2}, \eqref{decomposizione} and \eqref{gricia}
hypothesis \eqref{gricia*} holds, namely
\begin{equation}\label{griciabis}
\tetta_{-1}:=\frac{2^5 e K}{\a r s}\norma f^K\norma_{r,s}
\leq 1\,.
\end{equation}
Then
 there exists a real--analytic 
 symplectic change of coordinates
$$
\Psi_{-1}: D_{r_0}\times \T^n_{s_0} \ni (y^{(0)},x^{(0)})\
\to\ (y,x)\in D_r \times \T^n_s \,,\quad
r_0:=\frac34 r\,,\ \ \ s_0:=\left(1-\frac{1}{2K}\right)s\,,
$$
satisfying
\begin{equation}\label{argo2bis}
\noruno{y-y^{(0)}}
\leq 
\tetta_{-1}\frac{r}{16e}\,,
\qquad
\max_{1\leq i\leq n}|x_i-x^{(0)}_i|
\leq \
\tetta_{-1}\frac{s}{8 K}\,,
\end{equation}
such that 
\begin{equation}\label{olintobisbis}
H\circ\Psi_{-1}=:H_0=h(y^{(0)})+f_0(y^{(0)},x^{(0)})\,, \quad f_0=f^\flat+f_\varstar\,,\quad f^\flat:= P_\L f+\Tp\proiezione_\L^\perp  f\,,
\end{equation}
with
\begin{equation}\label{salaminabis}
\norma f_\varstar\norma _{r_0,s_0} 
\leq\ 
4\tetta_{-1} \norma f\norma _{r,s}
\,.
\end{equation} 
Recalling \eqref{decomposizione} and \eqref{olintobisbis}
we get
$$
f_0^K=f_\varstar^K
$$
and, by \eqref{salaminabis} and \eqref{griciabis},
\begin{equation}\label{lavatrice}
\norma f_0^K\norma _{r_0,s_0} 
\leq\ 
4\tetta_{-1} \norma f\norma _{r,s}
\leq
\frac{2^7 e K}{\a r s}\norma f\norma_{r,s}^2\,.
\end{equation}
Then, setting
\begin{equation}\label{corcira0}
\tetta_0:=\d \norma f_0^K\norma_{r_0,s_0}\qquad \mbox{ with}\qquad 
\d:=
\frac{2^5 e \, K^3}{\a r s}
\,,
\end{equation}
we have
\beq{cappuccino}
\tetta_0
\leq
\left( \frac{2^6 e  K^2}{\a r s}\norma f \norma_{r,s}\right)^2
\stackrel{\eqref{gricia}}\leq
(\tetta_\varstar/8)^2
\leq
\frac1{2^6}\ . 
\eeq
Finally, since $f_0^\flat-f^\flat=f_\varstar^\flat$ by \eqref{bemolle}
we get
\begin{equation}\label{salaminater}
\norma f_0^\flat-f^\flat
\norma _{r_0,s_0} 
\leq\ 
4\tetta_{-1} \norma f\norma _{r,s}
\stackrel{\eqref{griciabis}}
\leq 
\frac{2^7 e K}{\a r s}\norma f\norma_{r,s}^2
\stackrel{\eqref{gricia}}\leq
\frac{1}{4K}\tetta_\varstar\norma f\norma_{r,s}
\,.
\end{equation}

\nl
The idea is to construct $\Psi$ by applying $ K$ times Lemma~\ref{megara}. 

\nl
 Let 
\beqa{cornetto}
&&
\rho:=\frac{r}{4 K}\,, \qquad
\s:=\frac{s}{2K^2}
\,,
\nonumber\\
&& 
\nonumber\\
&&r_i:=\frac34 r- i \rho\,,\qquad
s_i:=\left(1-\frac{1}{2K}\right)s- i \s\,,\qquad
\norma \cdot\norma _i:=\norma \cdot\norma _{r_i,s_i}\,,
\eeqa
Fix $1\le j\leq K$ and
make the following {\bf inductive assumptions}: 

\nl
{\sl Assume that there exist, for $1\le i\le j$, real--analytic symplectic transformations
$$
\Psi_{i-1}:=X^1_{\phi_{i-1}} \ :\ D_{r_i}\times \T^n_{s_i}
\ni (y^{(i)},x^{(i)})\
\to\ (y^{(i-1)},x^{(i-1)})\in  
D_{r_{i-1}}\times \T^n_{s_{i-1}}\,,
$$
generated by a function $\phi_{i-1}=\phi^K_{i-1}$
with
\begin{equation}\label{maratonai}
\norma \phi_{i-1}\norma_{i-1}\leq \norma f^K_{i-1}\norma_{i-1}/\a\,,
\end{equation}
satisfying
\begin{equation}\label{argo3}
\noruno{y^{(i-1)}-y^{(i)}}
\leq \tetta_{i-1}\frac{r}{16eK}\,,\qquad
\max_{1\leq \ell\leq n}|x^{(i-1)}_\ell-x^{(i)}_\ell|\leq \tetta_{i-1}\frac{s}{8K^2}\,,
\end{equation}
such that
\begin{equation}\label{olintoj}
H_i:=H_{i-1}\circ \Psi_{i-1}
=:h+f_i=h+f_i^K+f_i^\flat
\end{equation}
satisfies, for $1\leq i\leq j$,  the estimates
\begin{equation}\label{pontina}
\tetta_i\leq \left(\frac{2^8 K^2\norma f\norma _{r,s}}{\a r s}\right)^{i+1} 
\stackrel{\eqref{gricia}}=\left(\frac{\tetta_\varstar}{8}\right)^{i+1} \,,
\qquad
\norma f_i^\flat-f^\flat_{i-1}\norma_i\le 
\frac{1}{\d}\left(\frac{\tetta_\varstar}{8}\right)^{i+1}
\,, 
\end{equation}
where}
\begin{equation}\label{corcira}
\tetta_i:=\d |f_i^K|_i
\ .
\end{equation}
Let us first show that  the inductive hypothesis is true for $j=1$
(which implies $i=1$). Indeed by \equ{cappuccino}
we see that  we
can apply Lemma~\ref{megara} with 
$f\rightsquigarrow f_0^K$ and $\check\tetta\rightsquigarrow\tetta_0=\d\norma f_0^K\norma_0$. Thus, we obtain the existence of $\Psi_0=X^1_{\phi_{{}_0}}$,
generated by a function $\phi_0=\phi^K_0$
with 
\begin{equation}\label{lavastoviglie}
\norma \phi_0\norma_{r_0,s_0}
\leq \frac{1}{\a}\norma f_0^K\norma_{r_0,s_0}
\stackrel{\eqref{lavatrice}}\leq
\frac{2^7 e K}{\a^2 r s}\norma f\norma_{r,s}^2
\,,
\end{equation}
satisfying \eqref{maratonai} 
and\footnote{Note also that $(f_0^K)^\flat=0$} \eqref{argo3},  so that 
$(h+f_0^K)\circ \Psi_0=:h+\tilde f_1$ and, by \equ{olintobis} and \equ{salamina},
\begin{equation}\label{recremisi}
 \norma \tilde f_1\norma_1 
 \le
  4  \tetta_0 \norma f_0^K\norma_0
  \stackrel{\eqref{cappuccino}}\leq \frac14\norma f_0^K\norma_0
   \stackrel{\eqref{lavatrice}}\leq
   \frac{2^5 e K}{\a r s}\norma f\norma_{r,s}^2
  \,.
  \end{equation}
  We have that $f_1= \tilde f_1+f_0^\flat\circ \Psi_0.$
  Then\footnote{Note that $(f_0^\flat)^K=0$ and $(f_0^\flat)^\flat=f_0^\flat.$}  
  \begin{equation}\label{anfizionia}
  f_1^K=\tilde f_1^K+(f_0^\flat\circ \Psi_0-f_0^\flat)^K\,,\qquad
  f_1^\flat-f_0^\flat=\tilde f_1^\flat+(f_0^\flat\circ \Psi_0-f_0^\flat)^\flat\,.
 \end{equation} 
  Write
  $$
  f_0^\flat\circ \Psi_0-f_0^\flat
  =
 ( f_0^\flat-f^\flat)\circ \Psi_0-( f_0^\flat-f^\flat)
  +
 ( f^\flat\circ \Psi_0-f^\flat-\{f^\flat,\phi_0\} )
  +\{f^\flat,\phi_0\}\,.
  $$
 By \eqref{delfi} (with
 $u\rightsquigarrow f_0^\flat-f^\flat,$
 $r \rightsquigarrow r_0,$ $s \rightsquigarrow s_0$) we have
 $$
 \norma
 ( f_0^\flat-f^\flat)\circ \Psi_0-( f_0^\flat-f^\flat)
 \norma_1
 \leq 
 2 \tetta_0 \norma f_0^\flat-f^\flat\norma_0
 \leq
\frac{2^4 e K}{\a r s} \norma f\norma _{r,s}^2
 $$ 
 by \eqref{cappuccino} and \eqref{salaminater}.
By \eqref{leuttra}
  with 
  $u \rightsquigarrow f^\flat,$
  $\phi\rightsquigarrow \phi_0,$
  $j\rightsquigarrow 1,$ 
  $\bar r \rightsquigarrow r_0,$
  $\bar s \rightsquigarrow s_0,$
  $$
 \norma  f^\flat\circ\Psi_0-f^\flat
-\{f^\flat,\phi_0\}\norma_1
\leq 
2\tetta_0
\norma
\{f^\flat,\phi_0\}\norma_{r_0-\rho/2,s_0-\s/2}
\leq
\frac{2^9 K^3}{\a^2 r^2 s^2}\norma f\norma_{r,s}^3 
\stackrel{\eqref{gricia}}\leq
\frac{ K}{4\a r s}\norma f\norma_{r,s}^2 
\,,
  $$
  by \eqref{cappuccino}, \eqref{lavastoviglie}
 and \eqref{tirinto} (with $f\rightsquigarrow \phi_0,$
 $g \rightsquigarrow f^\flat$).
 Analaogously by \eqref{tirinto} we get
$$
\norma \{f^\flat,\phi_0\}\norma_1
\leq
\frac{ 2^4 K^2}{e r s}\norma \phi_0\norma_0
\norma f\norma_{r,s}
\stackrel{\eqref{lavastoviglie}}\leq
\frac{ 2^{11} K^3}{\a^2 r^2 s^2}\norma f\norma_{r,s}^3
\stackrel{\eqref{gricia}}\leq
\frac{ K}{\a r s}\norma f\norma_{r,s}^2  \,.
$$  
Summarizing:
\begin{equation*}
\norma  f_0^\flat\circ \Psi_0-f_0^\flat\norma_1
\leq
\frac{ 2^6 K}{\a r s}\norma f\norma_{r,s}^2 \,.
\end{equation*}
 Then, by \eqref{recremisi} and \eqref{anfizionia} we get 
 \begin{equation}\label{frigo}
 \norma  f_1^K\norma_1\,,\ 
\norma  f_1^\flat-f_0^\flat\norma_1\ 
\leq\ 
\frac{ 2^7 K}{\a r s}\norma f\norma_{r,s}^2 
\end{equation}
checking \eqref{pontina} in the case $j=i=1.$
 
\nl
Now take $2\leq j\leq K$ and  assume that the inductive hypothesis holds true for 
$1\le i\le j$ and let us prove that it holds also for $i=j+1$.
By \eqref{pontina} and \eqref{gricia}
  we
can apply Lemma~\ref{megara} with 
$f\rightsquigarrow f_j^K$ and $\check\tetta\rightsquigarrow\tetta_j$. Thus, we obtain the existence of $\Psi_j=X^1_{\phi_{{}_j}}$,
generated by a function $\phi_j=\phi^K_j$
with 
\begin{equation}\label{lavastovigliej}
\norma \phi_j\norma_j
\stackrel{\eqref{maratonai}}\leq 
\frac{1}{\a}\norma f_j^K\norma_j
\stackrel{\eqref{corcira}}=
\frac{\tetta_j}{\a \d}
\,,
\end{equation}
 so that 
$(h+f_j^K)\circ \Psi_j=:h+\tilde f_{j+1}$ and, by \equ{olintobis} and \equ{salamina},
\begin{equation}\label{recremisij}
 \norma \tilde f_{j+1}\norma_{j+1} 
 \le
  4  \tetta_j \norma f_j^K\norma_j
  \stackrel{\eqref{corcira}}= \frac{4}{\d}\tetta_j^2
    \stackrel{\eqref{pontina}}\leq
    \frac{4}{\d}(\tetta_\varstar/8)^{2j+2}
    \stackrel{\eqref{gricia}}\leq
    \frac{1}{2^{3j-2}\d}(\tetta_\varstar/8)^{j+2}
    \leq
    \frac{1}{2^4\d}(\tetta_\varstar/8)^{j+2}
    \,,
  \end{equation}
  since $j\geq 2.$
  We have that $f_{j+1}= \tilde f_{j+1}+f_j^\flat\circ \Psi_j.$
  Then\footnote{Note that $(f_j^\flat)^K=0$ and $(f_j^\flat)^\flat=f_j^\flat.$}  
  \begin{equation}\label{anfizioniaj}
  f_{j+1}^K=\tilde f_{j+1}^K+(f_j^\flat\circ \Psi_j-f_j^\flat)^K\,,\qquad
  f_{j+1}^\flat-f_j^\flat=\tilde f_{j+1}^\flat+(f_j^\flat\circ \Psi_j-f_j^\flat)^\flat\,.
 \end{equation} 
Writing
$$
f_j^\flat=f^\flat + (f_0^\flat-f^\flat)+\sum_{h=1}^j f_h^\flat-f_{h-1}^\flat
$$
we have
\begin{eqnarray}
f_j^\flat\circ\Psi_j-f_j^\flat
&=&
\{ f^\flat, \phi_j \}
\nonumber
\\
&&+f^\flat\circ\Psi_j-f^\flat-\{ f^\flat, \phi_j \}
\nonumber
\\
&&
+(f_0^\flat-f^\flat)\circ\Psi_j-(f_0^\flat-f^\flat)
\nonumber
\\
&&
+
\sum_{i=1}^j 
\Big((f_i^\flat-f_{i-1}^\flat)
\circ\Psi_j -(f_i^\flat-f_{i-1}^\flat)\Big)
\label{nocchie}
\end{eqnarray}
where $\Psi_j=X^1_{\phi_j}$.
By \eqref{tirinto} with $f\rightsquigarrow \phi_j,$
 $g \rightsquigarrow f^\flat$,
 $r_0 \rightsquigarrow r_j$,
 $s_0 \rightsquigarrow s_j$, we get,
 by \eqref{maratonai} and \eqref{corcira},
$$
\norma \{f^\flat,\phi_j\}\norma_{j+1}
\leq
\frac{ 2^4 K^2}{e r s}
\norma \phi_j\norma_j
\norma f\norma_{r,s} 
\leq
\frac{ 2^4 K^2 \tetta_j}{e \a r s\d}\norma f\norma_{r,s} 
\stackrel{\eqref{gricia}}=
\frac{1}{e 2^4\d}(\tetta_\varstar/8)\tetta_j
\stackrel{\eqref{pontina}}\leq
\frac{1}{e 2^4\d}(\tetta_\varstar/8)^{j+2}
\,.
$$  
By \eqref{leuttra}
  with 
  $u \rightsquigarrow f^\flat,$
  $\phi\rightsquigarrow \phi_j,$
  $j\rightsquigarrow 1,$ 
  $\bar r \rightsquigarrow r_j,$
  $\bar s \rightsquigarrow s_j,$
  reasoning as above we get
  $$
 \norma  f^\flat\circ\Psi_j-f^\flat
-\{f^\flat,\phi_j\}\norma_{j+1}
\leq 
\tetta_j
\norma
\{f^\flat,\phi_j\}\norma_{r_j-\rho/2,s_j-\s/2}
\leq
\frac{\tetta_j}{4 e\d}(\tetta_\varstar/8)^{j+2}
\leq
\frac{1}{2^6 e\d}(\tetta_\varstar/8)^{j+2}
  $$
  by \eqref{pontina} and \eqref{gricia}.
By \eqref{delfi} (with
 $u\rightsquigarrow f_0^\flat-f^\flat,$
 $r \rightsquigarrow r_j,$ $s \rightsquigarrow s_j$) we have
 $$
 \norma
 ( f_0^\flat-f^\flat)\circ \Psi_j-( f_0^\flat-f^\flat)
 \norma_{j+1}
 \leq 
 2 \tetta_j \norma f_0^\flat-f^\flat\norma_j
 \leq
\frac{2^8 e K}{\a r s} \norma f\norma _{r,s}^2
\tetta_j
\leq 
\frac{1}{4\d}(\tetta_\varstar/8)^{j+2}
 $$ 
 by \eqref{salaminater}, \eqref{pontina},  \eqref{corcira0}
 and \eqref{gricia}.
Analogously, for $1\leq i\leq j,$ by \eqref{delfi} (now with
 $u\rightsquigarrow f_i^\flat-f_{i-1}^\flat$)
$$
\norma
(f_i^\flat-f_{i-1}^\flat)
\circ\Psi_j -(f_i^\flat-f_{i-1}^\flat)
\norma_{j+1}
 \leq 
 2 \tetta_j \norma f_i^\flat-f_{i-1}^\flat\norma_j
 \leq 
\frac{2}{\d}(\tetta_\varstar/8)^{j+i+2}
$$
 by  \eqref{pontina}.
 Then by \eqref{gricia}
 $$
 \norma \sum_{i=1}^j 
\Big((f_i^\flat-f_{i-1}^\flat)
\circ\Psi_j -(f_i^\flat-f_{i-1}^\flat)\Big)
\norma_{j+1}
\leq
\frac{2}{7\d}(\tetta_\varstar/8)^{j+2}\,.
 $$
Whence:
 $$
 \norma
 f_j^\flat\circ\Psi_j-f_j^\flat
 \norma_{j+1}
 \leq
\frac{4}{7\d}(\tetta_\varstar/8)^{j+2}\,.
 $$
 Then by \eqref{recremisij} 
 we get
 $$
\norma \tilde f_{j+1}\norma_{j+1}+\norma f_j^\flat\circ \Psi_j-f_j^\flat\norma_{j+1}
\leq
\frac{1}{\d}(\tetta_\varstar/8)^{j+2}\,.
 $$
By \eqref{anfizioniaj}
we get
\eqref{pontina} with $i=j+1.$
This completes the proof of the induction.

\nl
Now,  we can conclude the proof of Proposition \ref{pesce}. Set 
$$\Psi:=\Psi_{-1}\circ\Psi_0\circ\cdots\circ\Psi_{K-1}\ .
$$
Notice that, by \equ{cornetto},  $r_{K}=r/2=r_\varstar$ and $s_{K}=s(1-1/K)=s_\varstar$. By the induction, it is
\beq{bic}
H\circ \Psi=H_{{K}-1}\circ\Psi_{{K}-1}\stackrel{\equ{olintoj}_{K}}=h+f_{K}=:h+f^\flat+f_\varstar\,,
\eeq
with $f^\flat = \proiezione_\L f+ \Tp\proiezione_\L^\perp f$
(recall \eqref{senzanome}).
Note that by \eqref{pontina} and \eqref{cappuccino}
\begin{equation}\label{mantinea}
\sum_{i=1}^{K} \tetta_{i-1}
\leq
\sum_{i=1}^{K} (\tetta_\varstar /8)^i
\leq
\tetta_\varstar/7\,.
\end{equation}
Since 
$(y',x')=(y^{(K)},x^{(K)})$
by \eqref{argo2bis},
\eqref{argo3} and triangular inequality we get
\begin{eqnarray*}
\noruno{y'-y}
&\leq&
\noruno{y-y^{(0)}}+
 \sum_{i=1}^{K} 
 \noruno{y^{(i)}-y^{(i-1)}}
\leq
\frac{r\tetta_{-1}}{16 e }+
 \frac{r}{16 e K}\sum_{i=1}^{K} \tetta_{i-1}
 \\
&\stackrel{\eqref{mantinea}}\leq&
\frac{r}{16 e } \left(\tetta_{-1}+\frac{\tetta_\varstar}{7K}\right)
\stackrel{\eqref{griciabis}}\leq
\frac{r}{16 e } \left(\frac{\tetta_\varstar}{8K}+\frac{\tetta_\varstar}{7K}\right)
\,,
\end{eqnarray*}
then
 \eqref{argo0} follows
(the estimate on the angle being analogous).
\\
Since $T_K P_\L^\perp f^\flat=(f^\flat)^K=0$ (for any $f$, recall \eqref{decomposizione}) we have
\beqa{biro}
\norma T_K P_\L^\perp f_\varstar\norma _{r_\varstar,s_\varstar} 
&=&
 \norma f_{K}^K\norma _{K}\stackrel{\equ{corcira}}
 =
  \d^{-1} \tetta_{K}\stackrel{\equ{pontina}}
  \le
   \d^{-1} (\tetta_\varstar/8)^{{K}+1}
=
(\tetta_\varstar/8)^{{K}}\frac{8}{eK}\norma f\norma _{r,s}
\,,
\eeqa
proving the second estimates in \equ{pirati}. 

\nl
Finally, (using  that ${K}\ge 2$ and that $\tetta_\varstar \leq1$)
\beqano
\norma f_\varstar\norma _{r_\varstar,s_\varstar}&\stackrel{\equ{bic}}=&\norma f_{K}-f^\flat\norma _{K}\stackrel{\equ{decomposizione}}= 
\norma f_{K}^K+f_{K}^\flat-f^\flat\norma _{K}\le
\norma f_{K}^K\norma _{K}+\norma f_0^\flat-f^\flat\norma _0
+\sum_{i=1}^{K} \norma f_i^\flat-f_{i-1}^\flat\norma_i
\\
&\stackrel{\eqref{salaminater},\equ{pontina}}\le&
\norma f_{K}^K\norma _{K}
+\frac{1}{4K}\tetta_\varstar\norma f\norma_{r,s}
+\frac{1}{\d}\sum_{i=1}^{K}(\tetta_\varstar/8)^{i+1}
\\
&\stackrel{\equ{biro},\eqref{mantinea}}\le&
(\tetta_\varstar/8)^{{K}}\frac{8}{eK}\norma f\norma _{r,s}
+\frac{1}{4K}\tetta_\varstar\norma f\norma_{r,s}
+\frac{\tetta_\varstar^2}{56\d}
\le \frac{1}{K}\tetta_\varstar\norma f\norma_{r,s}\,,
\eeqano
which proves also the first estimate in \equ{pirati}.
\eproof

\section{Geometry of resonances}
\label{geometria}
We, first,  discuss the Covering Lemma in frequency space in a ball $B_M(0)\subset \{\o\in\real^n\}$
 and then we shall pull back through 
$\o^{-1}$ in the action domain. 
\\
We define a {\sl covering} $\{\O^i\}$ of $B_M(0)$ 
\beq{sanbernardo}
\O^0\cup\O^1\cup \O^2\supset B_M(0)\ ,
\eeq
as follows.

\begin{itemize}

\item[$\O^0$:]  The definition of the {\sl completely non--resonant zone $\O^0$} is nearly tautological:
\beq{neva0}
\O^0:=\{\o\in B_M(0):  \min_{k\in \genKO}|\o\cdot k|>\a/2  \}\ .
\eeq

\item[$\O^1$:] 
Recalling the definition of $\O^{1,k}$ in \eqref{sonno} we set
\begin{equation}\label{sonnosonno}
\O^1:=\bigcup_{k\in  \genKO} \O^{1,k}\,.
\end{equation}

\item[$\O^2$:]  The set $\O^2$ is the union of neighbourhoods of  exact double resonances\footnote{Recall \eqref{elrond}.} 
\beq{rkl}
R_{k,\ell}:=\{ \o\cdot k=\o\cdot \ell=0 \} \ ,\qquad k\in \genKO\ , \ell\in \genK,\ \ell\notin  \Z k\ ,
\eeq
namely:
\beq{palmer}
\O^2=\bigcup_{k\in \genKO} \bigcup_{ \ell\in \genK \atop \ell\notin  \Z k} \O^2_{k,\ell}
\eeq
where
\begin{equation}\label{wetton}
\O^2_{k,\ell}:=\{|\o\cdot k|<\a\}
\cap
\{|\proiezione_k^\perp \o|<M\}
\cap
\{|\proiezione_k^\perp \o \cdot \ell |\leq
3\a\K /|k|\}\,.
\end{equation}
\end{itemize}

\nl
Indeed, {\sl from these definitions, \equ{sanbernardo} follows immediately}. 

\giu
Next, let us point out the non--resonance properties satisfied by the frequencies in $\O^i$.

\begin{itemize}

\item[\rm (i)] If $0\neq \noruno{k}\le \KO$, then there exists $\bar k\in \genKO$ and a $0\neq j\in\integer$ such that 
$k=j \bar k$ and, therefore,
\beq{neva00}
\o\in\O^0 \quad \implies\quad 
|\o\cdot k|=|j||\o\cdot \bar k|\ge |\o\cdot \bar k|\ge \min_{k\in \genKO}|\o\cdot k|> \a/2\ .
\eeq

\item[\rm (ii)] Let $\o\in \O^{1,k}$ with $k\in\genKO$ and let $\ell\notin \Z k$, $\noruno{\ell}\le \K$.
Then,    there exist $j\in\Z\setminus\{ 0\}$
and $\ell'\in\genK$ such that $\ell=j\ell'$. Hence, 
\beqa{neva1}
|\o\cdot \ell|&=&|j|\, |\o\cdot \ell'| \ge |\o\cdot \ell'| =\Big| \frac{(\o\cdot k) (k\cdot \ell')}{|k|^2}+ \pk\o\cdot \ell'\Big|
\nonumber
\\
&\ge& |\pk\o\cdot \ell'|- \frac{\a \K}{|k|}> \frac{3 \a \K}{|k|}  - \frac{\a \K}{|k|} = \frac{2 \a \K}{|k|}\ .
\eeqa

\item[\rm (iii)] It remains to evaluate the measure of $\O^2$. To do this, we first prove the following

\begin{lemma}\label{lake}
 If $\o\in \O^2_{k,\ell}$ with 
 $k\in \genKO,\ \ell\in \genK$, $\ell \notin \Z k$, 
  then 
 \begin{equation}\label{melanzana}
 {\rm dist}(\o,R_{k,\ell})\leq 
 \sqrt{10}\, \a \K |k|\, |\ell| 
 \,.
\end{equation}
 Moreover,
\begin{equation}\label{zucchina}
 \meas (\O^2_{k,\ell}) \leq 
3\cdot 2^n \, M^{n-2} \a^2 \frac{\K}{|k|}\ .
\end{equation}
\end{lemma}
\proof
Let $v\in\R^n$ be the projection of $\o$ onto $R_{k,\ell}^\perp$, which is
the plane generated by $k$ and $\ell$
(recall that, by hypothesis, $k$ and $\ell$ are not parallel).
Then,   
\beq{sfinimento}
\dist (\o,R_{k,\ell})=\dist (v,R_{k,\ell})=|v|
\eeq
and
\begin{equation}\label{soldatino}
|v\cdot k|=|\o\cdot k|<\a\,, \qquad |\proiezione_k^\perp v \cdot \ell|
=|\proiezione_k^\perp \o \cdot \ell|
 \le 
3\a\K /|k|\,.
\end{equation}
Set
\beq{cacca}
h:=\pk \ell= \ell -\frac{\ell\cdot k}{|k|^2} k\,.
\eeq
Then, $v$ decomposes in a unique way as
$$v=a k+ b h$$
for suitable $a,b\in\R$.
By \eqref{soldatino},
\beq{goja}
|a|<\frac{\a}{|k|^2}\,,\qquad
|\pk v\cdot\ell|
=|bh\cdot \ell| \le 3\a\K /|k|\,,
\eeq
and 
$$
|h\cdot \ell|
\eqby{cacca}
\frac{ |\ell|^2 |k|^2-(\ell\cdot k)^2  }{|k|^2}
\ge \frac1{|k|^2}
$$
since $ |\ell|^2 |k|^2-(\ell\cdot k)^2$ is a positive integer (recall, that $k$ and $\ell$ are integer vectors not parallel).
Hence, 
\beq{velazquez}
|b|\le 3 \a \K |k|  \,,
\eeq
and \eqref{melanzana} follows since $|h|\le |\ell |$ and 
 $|v|=\sqrt{a^2 |k|^2+ b^2 |h|^2}\leq\sqrt{10} \a \K \, |k|\, |\ell |$.
\\
To estimate the measure of $\O^2_{k,\ell}$ we write $\o\in R_{k,\ell}$ as $\o=v+v^\perp$ with 
$v^\perp$ in the orthogonal
complement of the plane generated by $k$ and $\ell$. Since $|v^\perp |\le |\o|< M$  and $v$ lies in a rectangle of sizes of length $2\a/|k|^2$ and $6 \a \K |k|$ (compare \equ{goja} and \equ{velazquez})
we find
\beq{raimondi}
\meas(\O^2_{k,\ell})\le \frac{2\a}{|k|^2}\, (6 \a \K |k|)(2M)^{n-2}=3\cdot 2^n \, M^{n-2} \a^2 \frac{\K}{|k|}\ ,
\eeq
finishing the proof of Lemma~\ref{lake}.
\eproof
From \equ{palmer} and \equ{raimondi} it follows immediately (recall that $n\ge 2$) that
\beq{emerson}
\meas(\O^2)\le c M^{n-2} \a^2  \K^{n+1}\, \KO^{n-1}\ ,
\eeq
for a suitable constant $c$ depending only on $n$.

\end{itemize}

\proof {\bf (of Proposition \ref{beccalossi})}
Recalling \equ{neva0}, \eqref{sonnosonno} and \eqref{palmer},
for $i=0,1,2$ set
\beq{fischietto}
D^i:=\{y\in D: \o(y)\in \O^i\}
\ .
\eeq
Then \eqref{sanbernardo} implies \eqref{Di},
while
  \equ{neva00}, \equ{neva1} and \equ{emerson} imply immediately \equ{cipollotto}, \equ{cipollotto2} and\footnote{Recall the definition of $\bar L$ in Assumption {\bf A}, \S~\ref{mainresult}.} \equ{teheran4} respectively, proving Proposition~\ref{beccalossi}. \qed


\section{Averaging Theory}
\label{vento}

\noi
{\bf Assumption B} \\
{\sl 
Let $r,s>0$ and  let $h$ satisfy Assumption {\bf A} in \S~\ref{mainresult}.
\\
Let   $f:D_r\times \T^n_s\to\complex$  be a holomorphic function with
\begin{equation}\label{bada}
   \|f\|_{D,r,s}=1\,      
\end{equation}
and define
\beq{storia}
\ham(y,x):= h(y)+\e f(y, x)\ ,\qquad (y,x)\in D_r\times \T^n_s\ ,\  \e> 0\ . 
\eeq
Let
$\K$, $\KO$, $\nu$  and $\a$ be such that
\begin{equation}\label{islanda}
\K\ge  3 \KO\ge 6\ ,\qquad \quad
\nu\ge n+  2\,, \qquad \quad
\a:=\sqrt\e \K^\nu\,.
\end{equation}
For $k\in\genKO$, define}
\beqa{limone}
&\dst  \qquad r_{0}:=\frac{\a}{4L \KO}=\sqrt\e\frac{\K^\nu}{4L\KO}\ ; \qquad\quad & \dst r_k:=\frac{\a}{L |k|}=\sqrt\e\frac{\K^\nu}{L|k|}\ ,
\\
\label{tetta}
&\dst   \bar\tetta:= 
2^{14}n^{2n} \ \frac{L}{s^{2n+1}} \ \frac1{\K^{2\nu-2n-3}}
\ ;  \quad
&
\tetta
:=
2^{2n+10} n^{2n}  \frac{L}{s^{2n+1}}\frac{1}{\K^{2\nu-2n-3}}\ .
\eeqa


\noi
Putting together the Normal Form Lemma (Proposition~\ref{pesce}) and the Covering Lemma (Proposition~\ref{beccalossi}) there follows easily the following
{\sl averaging theorem for non--resonant and simply resonant zones}:

\begin{theorem}\label{raffo}
Let Assumption {\bf B} hold and assume that $\e$ satisfies
\eqref{dimenticavo} and
\beq{K*}
\K^{2\nu-n-4}\geq
2^{13+n}n^n \ \frac{L e^{s/2}}{s^{n+1}} 
\ ,
\eeq
then the following holds.

\nl
{\rm (i)}
There exists a symplectic change of variables
\begin{equation}\label{trota}
\Psi_{0}: D^0_{r_{0}/2}\times \T^n_{s(1-1/\KO)^2} \to 
D^0_{r_{0}} \times \T^n_{s(1-1/\KO)} 
\,,
\end{equation}
such that 
\begin{equation}\label{prurito}
\ham\circ\Psi_{0}
=h(y)+\e g^{{\rm o}}(y) +
\e f^{{\rm o}}_{\varstar\varstar}(y,x)\ ,  \qquad
\quad
\langle f^{{\rm o}}_{\varstar\varstar}\rangle=0\,,
\end{equation}
where $\langle \cdot \rangle=\proiezione_{\{0\}}$ denotes the 
average with respect to the angles $x$ and 
\beq{552}
\sup_{D^0_{r_{0}/2}}| g^{{\rm o}}-\langle f \rangle|
\leq
\bar\tetta\,, \qquad
\norma f^{{\rm o}}_{\varstar\varstar} \norma_{D^0,r_{0}/2,s(1-1/\KO)/2} 
\leq
2
\Big( \frac{2n\KO}{s}\Big)^n 
e^{-(\KO-3)s/2}
\,.
\eeq
{\rm (ii)}  $\dst D^1=\bigcup_{k\in\genKO} \Dk$ and for any $k\in \genKO$
there exists a symplectic change of variables
\begin{equation}\label{canarino}
\Psi_k: 
\Dk_{r_k/2}
\times \T^n_{s_{\varstar}} \to 
\Dk_{r_k} \times \T^n_{s(1-1/\K)} 
\,,
\qquad
s_{\varstar}:=s(1-1/\K)^2
\,,
\end{equation}
such that 
\begin{equation}\label{Hk}
\ham\circ\Psi_k
=h(y)+\e g^k(y,x) +
\e f^k_{\varstar\varstar} (y,x)
\end{equation}
where  
\begin{equation}\label{masiccio}
g^k=\proiezione_{{}_{k\Z}} g^k\,,\qquad
\proiezione_{{}_{k\Z}} f^k_{\varstar\varstar}=0\,,
\end{equation}
and
\begin{equation}\label{552bis}
\norma  g^k-\proiezione_{{}_{k\Z}} f\norma _{\Dk, r_k/2,s_{\varstar}}
\leq
 \tetta 
\,,
\qquad
\norma f^{k}_{\varstar\varstar} \norma _{\Dk,r_k/2,s(1-1/\K)/2} 
<
2
\Big( \frac{2n\K}{s}\Big)^n 
e^{-(\K-3)s/2}\ . 
\end{equation}
\end{theorem}

\rem\label{pianic}
(i) The functions $g^k$  and $\proiezione_{{}_{k\Z}}f$ depend, effectively,  only on one angle  $\sa\in\T^1$: more precisely, setting
\begin{equation}\label{sugna}
\left\{
\begin{array}l
F^k_j(y):=f_{jk}(y)\\ \ \\
G^k_j(y):=g^k_{jk}(y)
\end{array}
\right.
\qquad  \qquad
\left\{
\begin{array}l
\dst F^k(y,\sa):=\sum_{j\in\Z} F^k_j(y) e^{\ii j \sa}\, 
\\ \ \\ 
\dst G^k(y,\sa):=\sum_{j\in\Z} G^k_j(y) e^{\ii j \sa}
\end{array}
\right.
\end{equation}
we have (recall \eqref{dec})
\begin{equation}\label{palettabis}
(\proiezione_{{}_{k\Z}} f)(y,x)=
F^k(y,k\cdot x)\,,\qquad
g^k(y,x)=G^k(y,k\cdot x)\,.
\end{equation}
From  \eqref{552bis} and \eqref{stanco} it follows
\begin{equation}\label{cristina}
\norma  G^k-F^k\norma _{\Dk,r_k/2,|k|_{{}_1}s_{\varstar}} 
\leq  \tetta\ .
\end{equation}
The function $\sa\in\T^1_{\noruno{k}s_{\varstar}}\to G^k(y,\theta)$ will be called the {\bf effective potential} since, disregarding the small remainder $f^{k}_{\varstar\varstar}$, it governs the Hamiltonian evolution at simple resonances.

\nl
(ii) We have assumed that $\|f\|_{r,s}=1$ (see \equ{bada}), since this is the natural assumption in term of  genericity properties, however the Normal Form Lemma is formulated in term of the {\sl stronger} norm $\norma\cdot\norma$. We need therefore to restrict slightly the angle--analyticity domain in order to pass to the norm $\norma\cdot\norma$. This can be done through \eqref{battiato3}, which yields (for $r=r_0$ or $r=r_k$ and $K=\KO$ or $\K$)
\begin{equation}\label{bada2}
\norma f\norma_{r,s(1-1/K)}
\stackrel{\equ{battiato3},\equ{bada}}\le  
\left(\frac{2nK}{s}\right)^n
\,.
\end{equation}

\nl
(iii) The choice of $\a$ in \equ{islanda} is not restrictive (since it is done through the introduction of $\nu$, a new parameter) and it has the effect of making disappear $\e$ from the smallness conditions and from the definition of the smallness parameters $\bar\tetta$ and $\tetta$. \\
According to the choice of $\KO$ and $\K$  one will get different kind of statements. 
\erem

\rem
\label{irlanda} 
(i)
Observe that  $r_0\le r_k\le \sqrt{\e}\K^\nu/L$ so that assumption 
\equ{dimenticavo} ensures the necessary condition:
\beq{federico}
r_0\le r_k\le \frac{\sqrt{\e}\K^\nu}L\le r\ .
\eeq
(ii) The hypotheses of the Normal Form Lemma (Proposition~\ref{pesce}) concern {\sl a complex domain} $D_r$, while the non--resonance properties of the Covering Lemma (Proposition~\ref{beccalossi}) hold on {\sl real domains}.  
The following simple observation  allows to use directly the Covering Lemma:
  
\giu
{\sl If a set $D\subseteq \real^n$ is $(\a,K)$ non--resonant modulo $\L$ for $h$, then the complex domain $D_r$ is $(\a-LrK,K)$ non--resonant modulo $\L$, provided 
$L rK<\a$, where $L$ is  the Lipschitz constant of $\o$ on the complex domain $D_r$.}

\giu
Indeed, if $y\in D_r$ there exists $y_0\in D$ such that $|y-y_0|<r$ and 
$|\o(y_0)\cdot k|\ge \a$ for all $k\in\integer^n\bks\L$, $|k|_{{}_1}\le K$. Thus, for such $k$'s, one has
$$|\o(y)\cdot k|=|\o(y_0)\cdot k - (\o(y_0)-\o(y))\cdot k|\ge |\o(y_0)\cdot k|-L r K\ge \a-L rK\ . \qedeq
$$
\erem

\proof {\bf (of Theorem \ref{raffo})} 
(i): 
By Remark~\ref{irlanda}--(ii),
\eqref{cipollotto} and the choice of $r_0$ in \equ{limone}, 
the domain 
$D^0_{r_{0}}$ is $(\a/4,\KO)$ completely non--resonant (or non--resonant modulo the trivial lattice $\{0\}$) and, in view of \equ{bada2} and \equ{K*}, 
one can apply Proposition~\ref{pesce} to $\ham$ in
\eqref{storia}
with\footnote{Recall that the notation ``$ a \rightsquigarrow b$'' means ``with $a$ replaced by $b$''.}
\begin{eqnarray}
&&
f\rightsquigarrow\e f\ ,\quad 
D\rightsquigarrow D^0\ ,\quad
r\rightsquigarrow r_{0}\ ,\quad
\L\rightsquigarrow \{0\}\ , \quad
\a\rightsquigarrow \a/4\ ,
\nonumber
\\
&&
K\rightsquigarrow \KO\ ,\quad
s\rightsquigarrow s(1-1/\KO)\ .
\end{eqnarray}
Thus,  recalling \eqref{gricia}, using \equ{limone}  and that $\KO\ge 2$, one sees that 
\beqa{ELP}
\tetta_\varstar \rightsquigarrow \tetta_{0} &:=&
2^{15}\ 
\frac{L\KO^3  \norma f\norma_{r_0,s(1-1/\KO)}}{\K^{2\nu}  s(1-1/\KO)}
\nonumber\\
&\stackrel{\eqref{bada2},\equ{islanda}}<&
2^{16}\,  \frac{L \KO^3}{s\K^{2\nu} }\, \Big( \frac{2n\KO}{s}\Big)^n
\nonumber\\
&\stackrel{\equ{islanda}}\le& 
2^{13}n^n \ \frac{L}{s^{n+1}} \ \frac1{\K^{2\nu-n-3}}
\stackrel{\equ{K*}}\leq e^{-s/2}\le 1
\ ,
\eeqa
showing that  
 \eqref{gricia}   holds and also that $\bs\rightsquigarrow s(1-1/\KO)/2$ in \equ{shreck}.
Then,
by  \eqref{giovanni} and \equ{bada2}, one has:
\beqano
&&
\sup_{D^0_{r_{0}/2}}| g^{{\rm o}}-\langle f \rangle|
\leq
  \tetta _0\ \Big( \frac{2n\KO}{s}\Big)^n  
\stackrel{\equ{islanda}}\le 
 \Big( \frac{n\K}{s}\Big)^n  \tetta _0 
\stackrel{\eqref{ELP},\equ{tetta}}\leq
\bar\tetta
\,,
\nonumber
\\
&&
\norma f^{{\rm o}}_{\varstar\varstar} \norma_{D^0,r_{0}/2,s(1-1/\KO)/2} 
\leq
2  e^{-(\KO-2)s(1-1/\KO)/2}
\Big( \frac{2n\KO}{s}\Big)^n 
\leq 
2
\Big( \frac{2n\KO}{s}\Big)^n 
e^{-(\KO-3)s/2}
\,,
\eeqano
from which  \equ{552} follows.

\nl
(ii):
By Remark~\ref{irlanda}--(ii), the definition of $r_k$ in \equ{limone}
and \eqref{cipollotto2},
the domain 
$\Dk_{r_k}$ is 
$$(2\a \K/|k|-r_k L \K,\K)=(\a \K/|k|,\K)$$  
non--resonant modulo $\integer k$.
\\
Using again \equ{bada2}, 
we can apply Proposition~\ref{pesce}
with
\beqa{boston}
&& f\rightsquigarrow\e f \ , \quad
D\rightsquigarrow \Dk\ , \quad 
r\rightsquigarrow r_k\ , \quad 
\a\rightsquigarrow \a \K/|k|\ ,
\nonumber\\
&& K\rightsquigarrow \K
\ ,\quad
s\rightsquigarrow s(1-1/\K)
\ ,
\quad
\L\rightsquigarrow \Z k\ ,
\eeqa 
and (recall \eqref{gricia} and that $|k|\le \KO$)
\beqa{118}
\tetta_\varstar \rightsquigarrow \tetta_k 
&:=&
2^{11}\, \frac{L\K^2|k|^2 \e \norma f\norma_{r_k,s(1-1/\K)}}
{\a^2 s(1-1/\K)}
\nonumber\\ 
&\stackrel{\equ{islanda},\equ{bada2}}\leq& 
 2^{n+10} n^n  \frac{L}{s^{n+1}}\frac{1}{\K^{2\nu-n-4}}
\stackrel{\equ{K*}}
\leq e^{-s/2}
\le 1\ ,
\eeqa
showing that  
\eqref{gricia}  holds and also that $\bs\rightsquigarrow s(1-1/\K)/2$
in \equ{shreck}.
From  \eqref{giovanni}, \equ{118} and \equ{bada2} there follows \equ{552bis}; indeed 
\beqano
&&
\norma  g^k-\proiezione_{{}_{k\Z}} f\norma _{\Dk, r_k/2,s_{\varstar}}
\leq
\frac1{\K}
\Big( \frac{2n\K}{s}\Big)^n 
 \tetta_k 
\stackrel{\eqref{118},\equ{tetta}}\leq
\tetta
\,,
\nonumber
\\
&&
\norma f^{k}_{\varstar\varstar} \norma _{\Dk,r_k/2,s(1-1/\K)/2} 
\leq
2  e^{-(\K-2)s(1-1/\K)/2}
\Big( \frac{2n\K}{s}\Big)^n 
\leq 
2
\Big( \frac{2n\K}{s}\Big)^n 
e^{-(\K-3)s/2}
\,,
\eeqano
\qed

\section{Proofs of main results}
Under the above standing hypotheses, apart from a finite number of simple resonances 
{\sl the effective potential $G^k$ at simple resonances is close to a (shifted) cosine}:

\begin{proposition}\label{laura}
Let the assumptions of Theorem \ref{raffo} hold, let 
$k\in \genKO$ and 
let $G^k$ be as in \equ{sugna}, \equ{masiccio}. 
Then, if 
\beq{119}
\noruno{k}>3/s\ , 
\eeq
one has that
\beq{120}
\norma G^k - T_1 F^k\norma_{\Dk,r_k/2,2}\le  \tetta\, e^{s+5} e^{-\noruno{k}s} + 2^8 \, e^{-2\noruno{k}s}\ .
\eeq
\end{proposition}

\proof  Observe that by definition of\footnote{Recall \equ{labestia}.} $T_N$ and $T_N^\perp$,
\beq{allgreens} 
G^k-T_1 F^k= T_1G^k-T_1F^k+T_1^\perp G^k\ .
\eeq
Now, since $3/s<\noruno{k}\le \KO\le \K/3$,  
$$
\sup_{\Dk_{r_k/2}}
|G^k_{\!{}_{\pm 1}}-F^k_{\!{}_{\pm 1}}|\stackrel{\equ{cristina}}\le\tetta e^{-\noruno{k}s_{\varstar}} 
\stackrel{\equ{canarino}}\le \tetta e^{-\noruno{k}s (1-2/\K)} \le 
  \tetta e^s e^{-\noruno{k}s}
$$
so that
\beq{zorro}
\norma T_1G^k-T_1F^k\norma_{\Dk,r_k/2,2}=|G^k_1-F^k_1|e^2+ |G^k_{-1}-F^k_{-1}|e^2
<
2 e^2 e^s\tetta  e^{-\noruno{k}s}\ .
\eeq
Next, recalling \equ{sugna}, we have that
\beq{dondiego}
\norma T_1^\perp G^k\norma_{\Dk,r_k/2,2}
=\sum_{|j|\ge 2\atop j\in\integer} |g_{jk}| e^{2|j|}\le
\sum_{|j|\ge 2} |f_{jk}| e^{2|j|}+\sum_{|j|\ge 2} |g_{jk}-f_{jk}| e^{2|j|}\ .
\eeq
Let us estimate the two sums separately. 
Since $\|f\|_s=1$, $|f_\ell|\le e^{-\noruno{\ell}s}$ so that $|f_{jk}|\le e^{-|j|\noruno{k}s}$ and:
\beq{delavega}
\sum_{|j|\ge 2} |f_{jk}| e^{2|j|}\le
\sum_{|j|\ge 2}e^{-|j|\noruno{k}s} \ e^{2|j|} = 2\ \frac{e^{-2(\noruno{k}s-2)} }{1- e^{-(\noruno{k}s-2)}}
\le 4 e^4 e^{-2(\noruno{k}s)}\ ,
\eeq
where in the last inequality we used the assumption $\noruno{k}s>3> 2+\log 2$.\\
Then (again, because $\noruno{k}s>3$), we see that 
\beqa{bernardo} 
\sum_{|j|\ge 2} |g_{jk}-f_{jk}| e^{2|j|}
&=& \sum_{|j|\ge 2} |g_{jk}-f_{jk}| e^{|j|\noruno{k}s_\varstar} \ e^{-|j|\noruno{k}s_{\varstar}+2|j|}\nonumber\\
&\stackrel{\equ{cristina}}\le&\sup_{j\ge 2} \Big(e^{-j (\noruno{k}s_\varstar -2 ) }\Big)\  \tetta
\le e^{-2 (\noruno{k}s_\varstar -2 )} \ \tetta
\nonumber\\
&\le& 
 e^4 e^{-2\noruno{k} s (1-2/\K)} \tetta\le  \tetta\ e^{2s+4} e^{-2\noruno{k}s}\ .
\eeqa
Putting \equ{delavega} and \equ{bernardo} together, by \equ{allgreens} and \equ{119}, \equ{120} follows. \qed

\begin{proposition}\label{whisky}
Let the assumptions of Theorem \ref{raffo} hold; let $s>0,$ $0<\d\le 1$ and 
fix any  $0<\gamma\le 1$.
Assume 
\eqref{porticato2} and \eqref{porticato3}.
If  $k\in \genKO$ satisfies
\beq{porticato4}
\noruno{k}> 
 \ks{\d;\g}\ ,
\eeq
(with $\ks{\d;\g}$ defined in \equ{sidone})
then,
\beq{martha}
\norma G^k - T_1 F^k\norma_{\Dk,r_k/2,2}\le \g \, \d_k\ ,
\eeq
where
\beq{anatolia}
\d_k:=\d \noruno{k}^{-{n}} e^{-\noruno{k}s}\ .
\eeq
\end{proposition}

\rem\label{caffe!}
Conditions  \equ{porticato2} and \equ{porticato4} are stronger than the ones on $\nu$ in \equ{islanda}, \equ{K*} and  \equ{119}. In particular the assumptions of  Proposition~\ref{laura} hold. 
\erem 

\proof {\bf of Proposition \ref{whisky}}
As mentioned in the above remark, Proposition~\ref{laura} holds.
Let us estimate the two terms in \equ{120} separately. Recalling the definition of $\tetta$ in \equ{tetta} (and that $\noruno{k}\le \KO\le \K/3$), we find:
\beqa{ghevoi}
\tetta\, e^{s+5} e^{-\noruno{k}s}
 &\stackrel{\equ{tetta}}=& 
 e^{s+5} e^{-\noruno{k}s}
 2^{2n+10} n^{2n}  \frac{L}{s^{2n+1}}\frac{1}{\K^{2\nu-2n-3}}
 \nonumber
 \\
 &=&
  e^{s+5} 
 2^{2n+11} n^{2n}  \frac{L}{s^{2n+1}}\, \frac{\noruno{k}^{n}}{\K^{2\nu-2n-3}}\, 
\frac1{\g\d}\ \frac{\g\d_k}2\nonumber\\
&\le&
e^{s+5} 
 2^{n+11} n^{2n}   \frac{L}{s^{2n+1}}\, \frac{1}{\K^{2\nu-3n-3}}\frac1{\g\d}\ \frac{\g\d_k}2
\nonumber\\
&\stackrel{\equ{porticato2}}\le &
\frac{\g\d_k}2\,.
\eeqa
As for the second term in \equ{120}, we use the following calculus lemma, whose elementary check is left to the reader: 
\begin{lemma}\label{shhh}
If $a>2\log 2$, $0<\e<e^{-a^2/2}$ and $t>4\log \e^{-1}$, then $e^{-t}t^a<\e$.
\end{lemma}
Indeed, by the lemma (with $a\rightsquigarrow {n}$, $t\rightsquigarrow s\noruno{k}$ and $\e\rightsquigarrow s^{n} \g\d/2^9$)  and in view of  \equ{porticato3} and \equ{porticato4}, one has
\beqa{ghevoitu}
2^8 e^{-2\noruno{k}s}&=&(s\noruno{k})^{n} e^{-\noruno{k}s}\, \frac{2^9}{s^{n}\g\d}\cdot \frac{\g\d_k}2< \frac{\g\d_k}2\ .
\eeqa
The bounds \equ{ghevoi} and \equ{ghevoitu}  prove the claim. \qed

\noi
The  quantity $\d_k$ defined  in \equ{anatolia}  is a ``Fourier--measure''  for the non--degeneracy of analytic potentials $f$ holomorphic on  $\torus_s^n$, since such potential will have, in general, Fourier coefficients $f_k\sim e^{-s\noruno{k}}$.

\subsection{Positional potentials: proof of Theorem \ref{scozia}}

\nl
In order to conclude the proof of Theorem \ref{scozia} 
we need
the following

\begin{lemma}\label{edaje}
Let $s>0,$ $0<\d\le 1$ and 
fix any  $0<\gamma\le 1$.
Let the assumptions of Theorem~\ref{raffo} hold;  assume   \equ{porticato2}, \equ{porticato3} and that the positional potential
 $f\in  \cP_{s,{\tail_{\rm o}}}(\d)$ with the tail function $\tail_{\rm o}$
 defined in \eqref{porticato4}.
Then, for
$\ks{\d;\g}\le \noruno{k}\le \KO\leq \K/3$, one has
\beq{martina}
\sup_{y\in \Dk_{r_k/2}}\frac{\norma G^k(y,\cdot) - T_1 F^k(\cdot)\norma_{{}_2}}{|f_k|}\le \g 
\eeq
and
\begin{equation}\label{lothlorien}
\frac{1}{|f_k|}
\norma f^{k}_{\varstar\varstar} \norma _{\Dk,r_k/2,s(1-1/\K)/2} 
\leq
\frac{2^{10 n} n^{3n}}{ s^{3n}\d} 
e^{-\K s/8}
\ .
\end{equation}
\end{lemma}

\proof
Since Proposition \ref{whisky} holds, by \eqref{martha} and since
 $f\in  \cP_{s,{\tail_{\rm o}}}(\d)$ we get
\eqref{martina}.
By \eqref{552bis}
we get
\begin{eqnarray*}
&&\frac{1}{|f_k|}
\norma f^{k}_{\varstar\varstar} \norma _{\Dk,r_k/2,s(1-1/\K)/2} 
<
\frac{ \noruno{k}^{{n}} e^{\noruno{k}s}}{\d}
2
\Big( \frac{2n\K}{s}\Big)^n 
e^{-(\K-3)s/2}
\\
&&
\leq 
\frac{2^{n+1}n^n}{ s^n\d}\KO^{n} \K^n e^{-(\K-2\KO)s/2}
\stackrel{\eqref{islanda}}\leq
\frac{n^n}{ s^n\d} \K^{2n} e^{-\K s/6}
\,.
\end{eqnarray*}
Then
using that\footnote{Using that for $\a>0$ we have
$\max_{x>0} x^\a e^{-x}=(\a/e)^\a.$}
$$
\K^{2n} e^{-\K s/24}
\leq 
\left(\frac{48n}{s\, e} \right)^{2n}
\leq \frac{2^{10 n} n^{2n}}{s^n}\,,
$$
we prove \eqref{lothlorien}.
\eproof

\nl
Recalling the definition of $T_N$ given in \eqref{labestia},
we have that
$$T_1 F^k(\sa)=f_k e^{\ii \sa}+f_{-k} e^{-\ii \sa}= 2|f_k| \cos(\sa +\sa^{(k)})$$
for a suitable constant $\sa^{(k)}$.
Setting  
\begin{eqnarray*}
\mathtt G^k(y,\sa)
&:=&
\frac{G^k(y,\sa)-\cos(\sa+\sa^{(k)})}{2|f_k|}\,,
\\
\mathtt f^{k}(y,x)&:=&
\frac{f^{k}_{\varstar\varstar}(y,\sa)}{2|f_k|}\,,
\end{eqnarray*}
we get \eqref{HkTE}.
Finally Theorem \ref{scozia}  follows from Lemma \ref{edaje},
in particular \eqref{martinaTE} and \eqref{lothlorienTE} follow from
\eqref{martina} and \eqref{lothlorien}, respectively. \qed

\subsection{The general case ($y$-dependent potentials)}\label{denteazzurro}
For $k\in\Z^n\setminus\{0\}$ let $b_k>0$ such that
$$
\sum_{k\neq 0} b_k<\infty\,.
$$
For $\mathcal Z\subseteq \Z^n\setminus\{0\}$
we set 
$$
b_{\mathcal Z}:=\sum_{k\in \mathcal Z} b_k\,.
$$
For definiteness we will fix 
$$
b_k:= |k|^{-\frac{{n}}{2}}\,,
$$
but every other possible choice is fine.

\begin{proposition}\label{rivoluzione}
 Let $r,\mu>0$ and $\mathcal Z\subseteq \Z^n\setminus\{0\}.$
 For any $k\in\mathcal Z$ let
  $\gk (y)$ be holomorphic functions on the complex ball  $\{y\in\complex^n\, 
  \,:\, |y|<r\}$ with
 $$
 \sup_{|y|<r}|\gk (y)|\leq 1\,,\qquad
 \text{and}\qquad
 |\gk (0)|\geq \hat \delta_k>0 \,.
 $$
 Then, for every $y\in\real^n$ with $|y|<r/2e$, up, at most, to a  set of   
 measure\footnote{As usual $S^{n-1}:=\{y\in\R^n|\, \ |y|=1\}$.} 
 $$
\frac12 {\meas_{n-1} (S^{n-1})} b_{\mathcal Z}\left(\frac{r}{2e}\right)^{n}\mu
\,,
 $$
 we have
\begin{equation}\label{principe}
|\gk (y)|\geq \hat \delta_k \left(\frac{\mu b_k}{30 e^3}\right)^{\log1/\hat \delta_k}\,, \qquad
\forall\, k\in\mathcal Z\,.
\end{equation}
\end{proposition}
The proof relies on the following classical result in function theory (see, e.g., \cite{cartan}):
\begin{lemma} {\bf (Cartan's Estimate)}
 Assume that $\mathtt f:\C\to\C$ is holomorphic and bounded by $M>0$ on the 
 complex ball
 $|z|<2 e R$. If $|\mathtt f(0)|=1$  then, for $0<\eta<1$
 \begin{equation}\label{bottarga}
 |\mathtt f(z)|\geq \left(\frac{\eta}{15 e^3}\right)^{\log M}
\end{equation}
 for any $z\in\complex,$ $|z|<R$ 
 up to a set of balls of radii $r_j$ satisfying
 $$
 \sum_j r_j\leq \eta R\,.
 $$
\end{lemma}

\begin{remark}\label{piccolo}
 Note that 
 \eqref{bottarga}
holds in the complex ball  $|z|<R$ up to a set of measure smaller than
$ \pi \eta^2 R^2$.
Moreover it holds on the real interval $(-R,R)$ up to a set
of (real) measure $2\eta R.$
\end{remark}

\noi
\proof {\bf of Proposition~\ref{rivoluzione}}
Fix $k\in\mathcal Z.$
Fix $\xi=(\xi_1,\ldots,\xi_n)\in\R^n$ with $|\xi|=1$ and $\xi_1\geq 0.$
We apply Cartan's estimates simultaneously for every
$k\in\mathcal Z$ with
$$
\mathtt f(z)\rightsquigarrow \frac{\gk (z\xi)}{\gk (0)}\,,
\qquad
R \rightsquigarrow \frac{r}{2e}\,,
\qquad
M \rightsquigarrow \frac{1}{\hat \d_k}\,,
\qquad
\eta \rightsquigarrow \frac{\mu }{2}\,b_k.
$$
 By \eqref{bottarga}  estimate \eqref{principe} for the fixed $k$
 holds on the segment
 $\{y\xi \ :\ y\in(-r/2e,r/2e)\}$, up, at most, to a  set of   
 measure\footnote{Recall Remark \ref{piccolo}.}
 $$
  \frac{\mu }{2e}\, b_k r\ .
 $$
 Integrating on the half-sphere  $|\xi|=1$, $\xi_1\geq 0,$
 we get that 
 \eqref{principe} for the fixed $k$
 holds on the ball
 $|y|<r/2e$ up, at most, to a  set of   
 measure
 $$
 \frac12 {\meas}_{n-1}(S^{n-1})
 \left(\frac{r}{2e}\right)^{n}  \mu b_k\,.
 $$
 Summing on all $k\in\mathcal Z$
 we get that \eqref{principe} 
 holds for all $k\in\mathcal Z$.
\eproof

\medskip
\noi
Fix $0<\mu,\gamma<1.$
Define the following tail  function
\begin{eqnarray}\label{eowyn} 
&&\!\!\!\!\tail_*(\d;\g,\mu):=
\\
&&
\nonumber
\\
&&
\!\!\!\!
\frac{2^6n^2}{\ts}
\max\left\{
\ln^3 \frac{2^6{n}^2}{\ts}\,,\ 
\left(\ln\frac{30 e^3}{\d\mu}\right)\ln^2
\left(\frac{4}{\ts}\ln\frac{30 e^3}{\d\mu}\right)\,,\ 
\ln\frac{30 e^3}{\mu}\ln\frac{1}{\d}
\,,\ 
\ln\frac{2^{10}}{\d\g}
\right\}\,,
\nonumber
\end{eqnarray}
where
$$
\ts:=\min\{s,1\}\,.
$$
Fix  $y_0 \in D$
and assume that
\begin{equation}\label{galadriel}
f(y_0,\cdot)\in\cP_{s,\tail_*}(\d)\,.
\end{equation}
Set
\begin{equation}\label{rohan}
\gk(y):=f_k(y)e^{\noruno{k}s}\,.
\end{equation}
We have that
\begin{equation}\label{saruman}
\sup_{y\in\complex^n,\,|y-y_0|<r}|\gk(y)|\stackrel{\eqref{bada}}\leq 1\,,
\qquad
|\gk(y_0)|\geq \hat\d_k:=\d /\noruno{k}^{{n}}\,,
\qquad
\forall\, k\in\gen\,,\quad
|k|_{{}_1}> \tail_*(\delta)\,.
\end{equation}
Let $\mu>0$. Then by Proposition\footnote{With
$\gk(y)\rightsquigarrow \gk(y+y_0).$
} \ref{rivoluzione}
there exists a set\footnote{Both 
$ \mathcal D$ and $  B_{r/2e}(y_0)$ are real sets.}
 \begin{equation}\label{celeborn}
 \mathcal D\subseteq  B_{r/2e}(y_0)
 \quad {\rm satisfying}\quad
 \meas(B_{r/2e}(y_0)\setminus\mathcal D)\leq
\frac{b}{2} {\meas_{n-1} (S^{n-1})} \left(\frac{r}{2e}\right)^{n}\mu\,,
 \end{equation}
 with
 $$
 b:=\sum_{\noruno{k}>\tail_*(\d;\g,\mu)} \noruno{k}^{-{n}/2}
 \leq \sum_{k\neq 0} \noruno{k}^{-{n}/2}\,,
 $$
such that 
\begin{eqnarray}\label{gandalf}
&&|f_k(y)|e^{\noruno{k}s}=
|\gk(y)|\geq \d_k(\mu):=
\hat \d_k \left(\frac{\mu}{30 e^3 \noruno{k}^{{n}/2}}
\right)^{\ln 1/\hat \d_k}\,,
\nonumber
\\
&&\forall\, y\in \mathcal D\,, \ k\in\gen\,,
\quad
|k|_{{}_1}> \tail_*(\delta)\,.
\end{eqnarray}

\begin{theorem}\label{epitaph}
Let the assumption of Theorem \ref{raffo} hold.
Fix $0<\mu,\d<1/e^8$ 
and $0<\gamma<1$.
Assume that for some $y_0 \in D$
we have
$
f(y_0,\cdot)\in\cP_{s,\tail_*}(\d)$.
Set
\begin{equation}\label{fangorn}
\tilde \mu:=\mu/30 e^3\,,\quad
\tilde n:=2\nu-2n-3\,,\quad
\kappa:=2^{2n+10} n^{2n}  \frac{L}{s^{2n+1}}
\,.
\end{equation}
Assume that
\begin{equation}\label{frodo}
\K\geq 
\max\left\{
\KO^{\frac{2{n}^2}{\tilde n}\ln \KO}
\,,
\
\KO^{\frac{9}{\tilde n}\ln\frac{1}{\d\tilde \mu}}\,,
\
e^{\frac{4}{\tilde n}\ln\frac{1}{\d}\ln\frac{1}{\tilde\mu}}\,,
\
\left(\frac{4 e^{s+5}\kappa}{\d\g}\right)^{\frac{4}{\tilde n}}\,,\
\frac{2^5}{s}\ln^2\frac{1}{\d\mu}\,,\
\frac{2^{14}{n}^4}{s^2}
\right\}\,.
\end{equation}
If   $ k\in\genKO$
with $\noruno{k}> \tau_*(\d;\g,\mu)$
then
\begin{eqnarray}\label{feanor}
&&
\sup_{y\in(\Dk\cap\mathcal D)_{\hat r_k}}\,
\frac{\norma G^k(y,\cdot) - T_1 F^k(y,\cdot)\norma_{{}_2}}{|f_k(y)|}
\leq 
\gamma\,,
\\
&&
\sup_{y\in(\Dk\cap\mathcal D)_{\hat r_k}}\,
\frac{\norma f^{k}_{\varstar\varstar}(y,\cdot) \norma _{s(1-1/\K)/2}}{|f_k(y)|}
<
 \frac{4 e^{3s/2} n^n  }{\d s^n} 
  e^{-\K s/8}
\,,
\label{feanor2}
\end{eqnarray}
where  $ \mathcal D$ was defined in \eqref{celeborn} 
and\footnote{Note that by \eqref{federico} $\hat r_k\leq r/2.$} 
\begin{equation}\label{sauron}
\hat r_k:=\frac12 \min\{ r_k\,,\, \d_k(\mu)  \}\,.
\end{equation}
\end{theorem}

\proof
First we note that by \eqref{gandalf}, \eqref{saruman}, \eqref{sauron}
and Cauchy estimates
\begin{equation}\label{gandalf2}
|\gk(y)|\geq \frac12\d_k(\mu)\,,\qquad
\forall\, y\in \mathcal D_{\hat r_k}\,, \ k\in\gen\,,
\quad
|k|_{{}_1}> \tail_*(\delta)\,.
\end{equation}
By \eqref{120}, \eqref{rohan} and 	\eqref{gandalf2}
we have that for every $y\in \mathcal D\,, \ k\in\genKO$
\beq{120bis}
\frac{\norma G^k(y,\cdot) - T_1 F^k(y,\cdot)\norma_{{}_2}}{|f_k(y)|}
\leq 
\frac{2 e^{s+5} \, \tetta\ + 2^9 \, e^{-\noruno{k}s}}{\d_k(\mu)}\ .
\eeq
Then, in order to prove \eqref{feanor}, it is enough to show that
\begin{equation}\label{gondolin}
\frac{4 e^{s+5} \, \tetta}{\d_k(\mu)}\leq \gamma\,,\qquad
\frac{ 2^{10} \, e^{-\noruno{k}s}}{\d_k(\mu)}\leq \gamma\,.
\end{equation}
Let us consider the first inequality in \eqref{gondolin}.
Since  by \eqref{fangorn}  and recalling \eqref{tetta}
we have 
$\tetta=\kappa/\K^{\tilde n},$ then, recalling \eqref{saruman} 
and \eqref{gandalf},
for $\noruno{k}\leq \KO$
$$
\frac{4 e^{s+5} \, \tetta}{\d_k(\mu)}
\leq
\frac{4 e^{s+5} \kappa \KO^{{n}}}{\d\K^{\tilde n}}
\left(
\frac{\KO^{{n}/2}}{\tilde \mu}\right)^{\ln \KO^{{n}}\d^{-1}}
=
\frac{4 e^{s+5} \kappa}{\d} e^A\,,
$$
where
$$
A:=
\left(\frac{{n}}{2}\ln \KO+\ln\frac{1}{\tilde\mu} \right)
\left({n}\ln \KO+\ln\frac{1}{\d} \right)
+{n}\ln \KO
-\tilde n \ln \K
\,.
$$
Since
$$
A\leq
-\frac{\tilde n}{4}\ln\K
$$
by \eqref{frodo}, we obtain that
$$
\frac{4 e^{s+5} \, \tetta}{\d_k(\mu)}
\leq
\frac{4 e^{s+5} \kappa}{\d} e^{-\frac{\tilde n}{4}\ln\K}
\stackrel{\eqref{frodo}}\leq \g\,,
$$
proving the first estimate in \eqref{gondolin}.
\\
Regarding the second inequality in \eqref{gondolin}
we have
$$
\frac{ 2^{10} \, e^{-\noruno{k}s}}{\d_k(\mu)}
=\frac{2^{10}}{\d}e^B\,,
$$
with
$$
B:=
\left(\frac{{n}}{2}\ln \noruno{k}+\ln\frac{1}{\tilde\mu} \right)
\left({n}\ln \noruno{k}+\ln\frac{1}{\d} \right)
+{n}\ln \noruno{k}
-\noruno{k}s
\,.
$$
We note that 
$$
B\leq -\frac14\noruno{k}s
$$
by \eqref{eowyn}, 
indeed\footnote{
Note that $x/\ln x\geq \a$ if $x\geq \a \ln^2 \a$ and $\a\geq 5;$
analogously 
$x/\ln^2 x\geq \a$ if $x\geq \a \ln^3 \a$ and $\a\geq 2^6.$
}
$$
\frac{2}{s}{n}^2 \leq 
\frac{\noruno{k}}{\ln^2 \noruno{k}}\,,\qquad
\frac{4{n}}{s}\left(
\frac12\ln\frac{1}{\d}+\ln\frac{1}{\tilde\mu}+1
\right)
\leq
\frac{\noruno{k}}{\ln \noruno{k}}\,.
$$
Then
$$
  \frac{2^{10}}{\d}e^{-\frac14\noruno{k}s}
\stackrel{\eqref{eowyn}}\leq
\g\,.
$$
This proves \eqref{gondolin} and, therefore, 
completes the proof of \eqref{feanor}.
\\
Let us now show \eqref{feanor2}.
By \eqref{552bis}
and \eqref{gandalf2} we get
\begin{eqnarray*}
&&
\sup_{y\in(\Dk\cap\mathcal D)_{\hat r_k}}\,
\frac{\norma f^{k}_{\varstar\varstar}(y,\cdot) \norma _{s(1-1/\K)/2}}{|f_k(y)|}
<
\frac{2 e^{\noruno{k}s}}{\d_k(\mu)}
2
\Big( \frac{2n\K}{s}\Big)^n 
e^{-(\K-3)s/2}
\\
&&
\stackrel{\eqref{gandalf}}=
\frac{2 e^{\noruno{k}s}\noruno{k}^{n}}{\d} 
\left(\frac{30 e^3 \noruno{k}^{{n}/2}}{\mu}
\right)^{\ln (\noruno{k}^{n}/\d)}
2
\Big( \frac{2n\K}{s}\Big)^n 
e^{-(\K-3)s/2}
\\
&&
\stackrel{\eqref{islanda}}\leq
\frac{4 e^{3s/2} n^n }{\d s^n} 
\left(\frac{30 e^3 \K^{{n}/2}}{3^{{n}/2}\mu }
\right)^{\ln (\K^{n}/3^{n}\d)}
  \K^{2n}  e^{-\K s/6}
  \\
  &&
  =
  \frac{4 e^{3s/2} n^n  }{\d s^n} 
  e^{-\K s/8} e^{-Q}
\end{eqnarray*}
where
$$
Q:=\frac18 \K s
-\left({n}\ln\frac{\K}{3}+\ln\frac{1}{\d}
\right)
\left(\frac{{n}}{2}\ln\frac{\K}{3}
+\ln\frac{1}{\mu}+\ln 30 + 3
\right)
-2n\ln\K\,.
$$
Then \eqref{feanor2} follows if we prove that
$Q\geq 0.$
Recalling  \eqref{islanda}
we get\footnote{Using that $\ln^2 x\leq \sqrt x$ for $x\geq 2^{13}.$}
$$
Q\geq
\frac18 \K s
-8{n}^2\ln^2\K+2\ln^2\frac{1}{\d\mu}\geq 0
$$
by \eqref{frodo}.
\eproof

We rewrite Theorem \ref{epitaph}
in the fashion of Theorem \ref{scozia}.

\begin{theorem}\label{scoziaY}
Let $n\ge 2$, $0<s\le 1.$
Fix $0<\mu,\d<1/e^8$ 
and $0<\gamma<1$.
Consider a Hamilonian $H_\e(y,x)=h(y)+\e fy,(x)$ as in \eqref{ham}  
such that  $h$ satisfies the non--degeneracy Assumption~{\bf A} (\S~\ref{mainresult}) 
and  $f$ has norm one: 
$
   \|f\|_{D,r,s}=1.
   $
Assume that for some $y_0 \in D$
we have
$
f(y_0,\cdot)\in\cP_{s,\tail_*}(\d)$, with
$\tail_*=\tail_*(\d;\g,\mu)$
defined in \eqref{eowyn}.
Let
 $\K\ge  3 \KO\ge 6$ with $\KO$ 
satisfying 
\eqref{frodo} and \eqref{K*}.
Let  $\hat r_k$ as in \eqref{sauron}
and $ \mathcal D$ as in \eqref{celeborn} .
Finally assume that $\e$ satisfies 
\eqref{dimenticavo}.

\nl
Then,
for any   $k\in \genKO$ with
$\tail_*(\d;\g,\mu)\le \noruno{k}\le \KO$, 
there exists
a symplectic change of variables $\Psi_k$ as in \eqref{canarinoTE}
such that the following holds.

\nl
For every $y\in(\Dk\cap\mathcal D)_{\hat r_k}$
there exist a phase 
 $\sa^{(k)}(y)$ and 
 functions
 $\mathtt G^k(y,\cdot )\in\hol_2^1$
 and 
 $\mathtt f^k (y,\cdot )\in\hol_{s(1-1/\K)^2}^n $
 satisfying 
\begin{equation}\label{HkTEY}
\boxed{
\ham\circ\Psi_k
=:
h(y)+
2\e \,|f_k(y)|
\Big( \cos(k\cdot x +\sa^{(k)}(y))+
\mathtt G^k(y,k\cdot x)+
\mathtt f^k (y,x)
 \Big)
 }
\end{equation}
with 
\beq{martinaTEY}
\sup_{y\in(\Dk\cap\mathcal D)_{\hat r_k}}\,\norma \mathtt G^k(y,\cdot)\norma_2\le \g 
\eeq
and
\begin{equation}\label{lothlorienTEY}
\sup_{y\in(\Dk\cap\mathcal D)_{\hat r_k}}\,
\norma \mathtt f^{k} \norma _{s(1-1/\K)/2} 
\leq
 \frac{4 e^{3s/2} n^n  }{\d s^n} 
  e^{-\K s/8}
\ .
\end{equation}
\end{theorem}

\proof
It directly follows from Theorem \ref{epitaph}.
We only note that $\sa^{(k)}(y)$
is defined such that
$$
|f_k(y)|
 \cos(k\cdot x +\sa^{(k)}(y))=T_1 F^k(y,x)\,,
$$
while
\beqano
\mathtt G^k(y,x) &:=&
\frac{G^k(y,\cdot) - T_1 F^k(y,\cdot)}{|f_k(y)|}\,,
\\
\mathtt f^k(y,x) &:=&
\frac{f^{k}_{\varstar\varstar}(y,x) }{|f_k(y)|}\,,
\eeqano
Note that $\sa^{(k)}(y),$ $\mathtt G^k(y,x)$ and $\mathtt f^k(y,x)$
are not analytic in $y$ (due to the presence of $|f_k(y)|$), but, obviously, $\ham\circ\Psi_k$ is real--analytic in $x$ {\sl and} $y$.
\eproof

\giu
{\bf Acknowledgments}
We are grateful to  V. Kaloshin and A. Sorrentino for useful discussions.

\appendix

\section{An elementary result in linear algebra}\label{appendicedialfonso}
\begin{lemma}\label{alfonso}
Given $k\in\Z^n,$ $k\neq 0$
there exists a matrix $A=(A_{ij})_{1\leq i,j\leq n}$ with integer entries such that
$A_{nj}=k_j$ $\,\forall\ \, 1\leq j\leq n$,
$\det A=d:={\rm gcd}(k_1,\ldots,k_n)$, and 
$|A|_{{}_\infty}= |k|_{{}_\infty}$.
\end{lemma}

\proof
The argument is by induction over $n$.
For $n=1$ the lemma is obviously true. For $n=2$, it follows at once from\footnote{The first statement in this formulation of Bezout's Lemma is well known 
and it can be found in any textbook on elementary number theory; the estimates on $x$ and $y$ are easily deduced from the well known fact that 
given a solution $x_0$ and $y_0$ of the equation $a x + b y =d$, all other solutions have the form $x= x_0+k (b/d)$ and $y=y_0- k (a/d)$ with $k\in \integer$
and by choosing $k$ so as to minimize $|x|$.
} 

\nl
{\bf Bezout's Lemma} {\sl Given two integers $a$ and $b$ not both zero, there exist two integers $x$ and $y$ such that $a x+ b y = d:={\rm gcd}(a,b)$, and such that 
$\max\{|x|,|y|\}\le \max\{|a|/d, |b|/d\}$.}

\nl
Indeed, if $x$ and $y$ are as in Bezout's Lemma with $a= k_1$ and $b=k_2$ one can take $A= \begin{pmatrix}y & -x\\ k_1 & k_2 \end{pmatrix}$.
Now, assume, by induction for $n\ge 3$ that the claim holds true for $(n-1)$ and let us prove it for $n$.   
Let $\bar k = (k_1,...,k_{n-1})$ and $\bar d={\rm gcd}(k_1,...,k_{n-1})$ and notice that ${\rm gcd}(\bar d, k_n)=d$. By the inductive assumption, there exists a matrix $\bar A=\begin{pmatrix} \tilde A\\ \bar k \end{pmatrix}\in
 {\rm Mat}_{(n-1)\times (n-1)}(\Z)$  with  $\tilde  A\in {\rm Mat}_{(n-2)\times (n-1)}(\Z)$, such that $\det \bar A=\bar d$ and $|\bar A|_{{}_\infty}=|\bar k|_{{}_\infty}$.
Now,  let $x$ and $y$ be as in Bezout's Lemma with $a= \bar d$, and $b=k_n$.
We claim that $A$ can be defined as follows:
\beq{giaccherini}
A=\begin{pmatrix} & \tilde k & \ & \tilde x \\ & \bar A &\  & \begin{pmatrix} 0 \\ \vdots \\ 0 \\k_n\end{pmatrix} \end{pmatrix}\ ,\qquad \tilde k= (-1)^n y\,  \frac{\bar k}{\bar d}\ ,\qquad \tilde x:= (-1)^{n+1} x \ .
\eeq
First, observe that since $\bar d$ divides $k_j$ for $j\le (n-1)$, $\tilde k \in \integer^{n-1}$. Then, expanding the determinant of $A$ from last column, we get
\beqano
\det A& =& (-1)^{n+1}\tilde x \det \bar A + k_n \det \begin{pmatrix} \tilde k \\ \tilde A\end{pmatrix}\\ 
& =& 
(-1)^{n+1}\tilde x \, \bar d + k_n (-1)^{n-2} \det  \begin{pmatrix} \tilde A \\ \tilde k\end{pmatrix} \\
&=& (-1)^{n+1}\tilde x\, \bar d + k_n  (-1)^{n-2}  (-1)^n \frac{y}{\bar d} \det \bar A\\
&=& 
x \bar d + k_n y = d\ .
\eeqano
Finally, by Bezout's Lemma, we have that $\max\{|x|,|y|\}\le \max\{\bar d/d , |k_n|/d\}$, so that
$$
|\tilde k|_{{}_\infty} = |y| \frac{|\bar k|_{{}_\infty} }{\bar d} \le \frac{|\bar k|_{{}_\infty}}{d}\le |k|_{{}_\infty}\ ,\quad
|\tilde x|= |x|\le \frac{|k_n|}{d}\le |k|_{{}_\infty}\ ,
$$
which, together with $|\bar A|_{{}_\infty}=|\bar k|_{{}_\infty}$,  shows that $|A|_{{}_\infty}=|k|_{{}_\infty}$. \qed

\comment{
\section{The mechanical system case}
\nl
Near an exact  simple resonance $\{\o\cdot k=0\}$, with $k\in\gen$,  the averaged Hamiltonian will be the projection of
the perturbation $f$ on the one--dimensional Fourier lattice generated by $k$, $k\Z$. Such projection is exactly the 1d--Fourier projection $\pi_{k\Z}$ of $f(y,\cdot)$ defined in \equ{pale} evaluated at the ``angle'' $\sa=x\cdot k$:
\beq{paley}
\proiezione_{k\Z}f(y,x)=F^k(k\cdot x,y)\ ,\qquad F^k(\sa,y):=\sum_{j\in \integer\backslash\{0\}} f_{jk}(y)e^{\ii j \sa}\ .
\eeq
Thus, the new ``natural'' angles  after averaging, will involve a new special angle, say $\psi_n=k\cdot x\in\torus^1$, which will be the effective secular angle on which depend the averaged Hamiltonian (disregarding the remainder, which, of course, will depend on other $n-1$ angles).  There are, in general, many linear symplectomorphisms  on $\R^n\times \torus^n$ so that  $x\in\torus^n \mapsto \psi\in\torus^n $ with $\psi_n=x\cdot k$.  Indeed, if $\Ak \in SL(n,\Z)$  and has one raw given by\footnote{Here, $k$ is  a row vector (normally we do not distinguish between row and column vectors since it will be clear from context); the further specification on the norm of $\hAk$, is also, always achievable, compare, e.g., Lemma~\ref{alfonso} in Appendix~\ref{appendicedialfonso}.} $k$, i.e.,  
\begin{equation}\label{scimmia}
\left\{\begin{array}l
\dst \Ak =\binom{\hAk }{k}\in\ {\rm SL}(n,\Z)\ ,
\\ \ \\ 
\hAk \in{\rm Mat}_{(n-1)\times n}(\Z)\ ,
\end{array}\right.
\qquad\quad 
\big(|\hAk |_{{}_\infty}\leq |k|_{{}_\infty}\big)\ ,
\end{equation}
then the map 
\beq{fulmini}
(y,x)\in\real^n\times \torus^n\mapsto (J,\psi)\in \real^n\times \torus^n\ :\qquad 
\left\{
\begin{array}l
\psi= \Ak x\ ,\\ \ \\ J=\Ak^{-T}y\ ,
\end{array}\right.
\eeq
which is such that $\psi_n=k\cdot x$, 
is a global symplectomorphism.
\\
The description of the covering is more conveniently described in terms of the new symplectic variables $(J,\psi)$. 

\beq{chemaroni}
G^k(y,\sa)=2|f_k(y)|\big(-\cos(\sa +\sa^{(k)}(y))+ R^k(y,\sa)\big)
\eeq
one with
where $\sa^{(k)}=\pi + \arctan(\Im f_k/\Re f_k)$ and $R^k$ (defined by such relation) is
satisfies the bound
\beq{terame}
\norma R^k\norma _{\tDk,r_k/2,2}\le \frac14\ . 
\eeq
}


\begin{thebibliography}{99}

\footnotesize

\bibitem{A} V.~I. Arnold.
Instability of dynamical systems with several degrees of freedom. Dokl. Akad. Nauk SSSR 156, 9--12 (1964). Engl. transl.: Sov. Math., Dokl. 5, 581--585 (1964)

\bibitem{AKN}
V.~I. Arnold, V.~V.  Kozlov, and A.~I. Neishtadt.
Mathematical aspects of classical and celestial mechanics,\emph{
volume~3 of Encyclopaedia of Mathematical Sciences}.
Springer-Verlag, Berlin, third edition, 2006.
[Dynamical systems. III], Translated from the Russian original by E. Khukhro.



\bibitem{BC} L. Biasco, and L. Chierchia. KAM Theory for secondary tori, 
 arXiv:1702.06480v1 [math.DS] 

\bibitem{BClin} L. Biasco, and L. Chierchia. 
On the measure of Lagrangian invariant  tori in nearly--integrable mechanical systems. 
Rend. Lincei Mat. Appl. {\bf 26} (2015), 1--10 




%

%






\bibitem{delshamsG} 
A. Delshams; P. Guti\'errez,
Effective stability and KAM theory.  
J. Differential Equations 128 (1996), no. 2, 415--490



\bibitem{giorgilli}
A. Giorgilli; L. Galgani,
Rigorous estimates for the series expansions of Hamiltonian perturbation theory. 
Celestial Mech. 37 (1985), no. 2, 95--112


\bibitem{GCB} M. Guzzo, L. Chierchia and G. Benettin. The steep Nekhoroshev Theorem, 
Commun. Math. Phys. 342, 569-601 (2016) 

\bibitem{HK}
B. R. Hunt, V. Y. Kaloshin, 
Prevalence, chapter 2, Handbook in dynamical systems, edited by H. Broer, F. Takens, B. Hasselblatt, Vol. 3, 2010, pg. 43-87

\bibitem{L} V. F. Lazutkin, Concerning a theorem of Moser on invariant curves (Russian), Vopr. Dinamich. Teor. Rasprostr. Seism. Voln. 14 (1974), 109-120.





\bibitem{cartan} B. Y. Levin. Lectures on Entire Functions, 
American Mathematical Soc., 1996.

%
%
%

\bibitem{MNT} A.G. Medvedev, A.I.  Neishtadt, D.V. Treschev, 
Lagrangian tori near resonances of near--integrable Hamiltonian systems, 
Nonlinearity, {\bf 28}:7 (2015), 2105--2130


\bibitem{Nei}
A. I. Neishtadt, Estimates in the Kolmogorov theorem on conservation of condition- ally periodic motions, J. Appl. Math. Mech. 45 (1981), no. 6, 766-772




\bibitem{nek}
N.N. Nekhoroshev. An exponential estimate of the time of stability of nearly- integrable Hamiltonian systems I, Uspekhi Mat. Nauk 32 (1977), 5--66; Russian Math. Surveys 32 (1977), 1--65
%
%
%
%



\bibitem{poschel1982} 
J. P\"oschel, Integrability of Hamiltonian systems on Cantor sets, Comm. Pure Appl. Math., v. 35 (1982), no. 1, 653-695


\bibitem{poschel}
J. P\"oschel,  
{Nekhoroshev estimates for quasi--convex Hamiltonian systems.} 
Math. Z. {\bf 213}, pag. 187 (1993).





\end{thebibliography}
\end{document}